\documentclass{article}
\usepackage{amssymb,amsmath}
\usepackage{anysize}
\usepackage{psfrag}
\newtheorem{thm}{Thm.}[section]
\newtheorem{def.}{Def.}[section]

\newtheorem{prop}{Prop.}[section]
\newtheorem{cor}{Cor.}[section]

\newtheorem{ex}{Example}[section]
\numberwithin{table}{section}

\begin{document}
\title{On Finite Racks and Quandles}
        \author{Pedro Lopes\\
        Departamento de Matem\'atica\\
            Instituto Superior T\'ecnico\\
            Av. Rovisco Pais\\
        1049-001 Lisboa\\
            Portugal\\
        \texttt{pelopes@math.ist.utl.pt}\\
\and
        Dennis Roseman\\
        Mathematics Department\\
        The University of Iowa\\
        MacLean Hall\\
        Iowa City, IA 52242\\
        USA\\
        \texttt{roseman@math.uiowa.edu}
}
\date{August 26, 2004}
\maketitle

\begin{abstract}
We revisit finite racks and quandles using a perspective based on
permutations which can aid in the understanding of the structure.
As a consequence we recover old results and prove new ones. We
also present and analyze several examples.
\end{abstract}

\section{Introduction} \label{S:intro}

The algebraic structure known as quandle was introduced
independently by Joyce and Matveev (\cite{dJoyce} and
\cite{sMatveev}) in order to obtain invariants of knotted circles
i.e., invariants of embeddings of the circle in three space
(\cite{lhKauffman}). Given any knotted circle, the idea is to
associate a quandle to it - the knot quandle. As was later
remarked by Fenn and Rourke (\cite{rFenn0}) this procedure
generalizes to higher dimensions, in particular to embeddings of
closed surfaces in four space - the so-called knotted surfaces
(\cite{dRoseman1}, \cite{CetS}). Given the similarities of these
two procedures we will next describe how to associate the knot
quandle to the embedding in both cases i.e., embeddings of circles
in three space and embeddings of closed surfaces in four space.
Consider the embedding under study. Upon convenient projection of
it on a plane ($3$-hyperplane, resp.), the crossing points (lines
of self-intersection, resp.) split the projection into arcs
(regions, resp.).

A presentation of the knot quandle is obtained by regarding the
arcs (regions, resp.) as generators and by reading specific
relations at the crossings (lines of self-intersection, resp.).
Thanks to the defining axioms of quandles this association is
essentially insensitive to the Reidemeister moves (Roseman moves,
resp.) thus making the knot quandle an invariant for knotted
circles (for knotted surfaces, resp.). Moreover, in the case of
knots, the knot quandle is a classifying invariant modulo
orientation of the ambient space. Given that the knot quandle is
an invariant modulo orientation of the ambient space, we will take
a knot to be a representative of each of these equivalence
classes; in particular, the trefoil knot and its mirror image
correspond to the same knot.

As in the case of group presentations, these quandle presentations
by themselves are difficult to work with. On the other hand, the
na\"{\i}ve approach of counting homomorphisms from the knot
quandle to specified target quandles (counting colorings of the
knot quandle by specified labelling quandles, \cite{rFenn0},
\cite{jsCetal}, \cite{pLopes}) has proved to be extremely
effective. As a matter of fact, with only ten labelling quandles
all but less than $3 \%$ of the pairs of prime knots of up to ten
crossings are told apart by counting colorings (\cite{DL}). A
relevant success is also obtained when trying to tell apart
knotted surfaces via counting colorings (\cite{BL}).

There is a related algebraic structure known as rack of which the
quandles are particular cases. Racks are invariants of framed
knots (\cite{rFenn0}). Moreover, racks are invariants of knotted
surfaces (\cite{dRoseman2}).

Notwithstanding the success of these applications of quandles
and/or racks to knot theory, there are simple questions that
remain unanswered. Here are some of them.  Are there knots which
are not distinguished by counting colorings? If any two knots are
distinguished by counting colorings, given an integer $c>2$, find
an integer $n\sb c>2$ such that all prime knots of crossing number
less than or equal to $c$ are distinguished via counting colorings
with labelling quandles of order less than or equal to $n\sb c$.
Can the class of labelling quandles be further restricted to
finite Alexander quandles? Or even to finite linear Alexander
quandles?

We believe that in order to answer these and other related
questions more has to be known about the structure of racks and
quandles. In this paper we focus on racks and quandles as
sequences of permutations. In Section \ref{S:prel} we describe the
passage into regarding racks and quandles as sequences of
permutations and give basic examples of racks from this point of
view. We rewrite the axioms and prove easily basic properties in
Section \ref{S:perm}. In Sections \ref{S:inv1} and \ref{S:inv2} we
use elementary facts about permutations to derive invariants of
racks and quandles. In Section \ref{S:indecomptrans} we briefly
discuss indecomposability and transitivity of racks and quandles.
In Section \ref{S:orbs} we introduce $\Delta$-orbits which will be
useful in the proofs in Section \ref{S:racksfromprof} where we
identify the racks that realize given profiles. Throughout the
text we present and discuss several examples of racks and
quandles.

\section{Preliminaries} \label{S:prel}

A {\bf rack} (\cite{dJoyce}, \cite{sMatveev}, \cite{Bries},
\cite{lhKauffman}, \cite{rFenn0}, \cite{pD}, \cite{AG}, for
instance) is an algebraic structure whose binary operation,
denoted $\ast $, is right-invertible and self-distributive.
Denoting $X$ the underlying set, right-invertibility states that
for any $i$, $j$ in $X$, there is a unique $x$ in $X$ satisfying
the equation $x\ast j = i$. Hence, for any $j$ in $X$, the
assignment
\[
\begin{array}{clcr}
    \mu\sb j :& X \longrightarrow X\\
              & x \longmapsto x\ast j
\end{array}
\]
is a permutation of $X$; it is surjective because there is a
solution to each equation referred to above and it is injective
because that solution is unique, see Section \ref{S:perm} for
details. We will then regard a rack of order $n$ as a finite
sequence of $n$ permutations from $S\sb n$ (the group of
permutations of $n$ objects) one for each element of $X$, in the
sense just discussed. In particular, the $j$-th column of the
multiplication table of the rack under study is the image of the
permutation corresponding to $j$, for each $j$ in $X$. Since there
is a second axiom in the definition of rack, it turns out that not
any sequence of $n$ permutations from $S\sb n$ gives rise to a
rack. In fact, this second axiom, self-distributivity, amounts to
saying that the $\mu\sb{j}$'s preserve the $\ast$ operation i.e.,
$\mu\sb{j}(i\ast i')=\mu\sb{j}(i)\ast \mu\sb{j}(i')$ and since the
$\mu\sb{j}$'s are bijections then this amounts to saying that each
$\mu\sb{j}$ defines an automorphism on the rack under study. We
will next illustrate these ideas with examples. Our main reference
for notation and facts about permutations will be \cite{jjRotman}.

\subsection{Preliminary Examples} \label{S:prelex}

\begin{ex} For integer $n$ greater than $2$, the {\bf dihedral rack} of
order $n$, denoted $R\sb n$, is the set formed by the first $n$
positive integers, $\{1,2,\dots ,n\}$, equipped with the operation
$i\ast j := 2j-i$, mod $n$.\end{ex}

The multiplication table for $R\sb 3$ is displayed in Table
\ref{Ta:R3}.
\begin{table}[h!]
\begin{center}
    \begin{tabular}{| c || c | c | c |}\hline
    $R\sb 3$  &   1     &      2  &     3  \\ \hline \hline
            1 &   1     &      3  &     2   \\ \hline
            2 &   3     &      2  &     1   \\ \hline
            3 &   2     &      1  &     3   \\ \hline
    \end{tabular}
\caption{$R\sb 3$ Multiplication Table}\label{Ta:R3}
\end{center}
\end{table}
The permutation corresponding to $1$ is then
\[
\mu\sb 1 =
\begin{pmatrix}
1 & 2 & 3\\
1 & 3 & 2
\end{pmatrix} = (1)(2 3)
\]
where $(1)(2 3)$ is the permutation expressed as products of
disjoint cycles, which is the way permutations will be written,
henceforth; note also that $()$ denotes the identity permutation.
Analogously, $\mu\sb 2 = (1 3)$, and $\mu\sb 3 = (1 2)$, omitting
$1$-cycles. The sequence of permutations of $R\sb 3$ is then
$\big( (2 3), (1 3), (1 2)\big)$. In this way, $1\ast 1 = \mu\sb
1(1)=1$, $1\ast 2 = \mu\sb 2(1)=3$, $\mu\sb j(i)=i\ast j=2j-i$,
mod $3$, for any $i$, and $j\in \{1, 2, 3\}$ thereby recovering
the rack operation on $\{1, 2, 3\}$. For general $n$, the $i$-th
term of the sequence of permutations of $R\sb{n}$ is, for any $i$
in $\{1, 2,\dots ,n\}$, the product of transpositions,
\[
\mu\sb{i}=\prod\sb{j=1}\sp{[\frac{n-1}{2}]} (\; i+j \quad i-j\; ) \qquad \texttt{mod n}
\]

\begin{ex} $T\sb n$ is the {\bf trivial rack} of order $n$. The
underlying set is the set of the $n$ first positive integers $\{1,
2, \dots ,n\}$ and $i\ast j := i$ for any $i$, $j\in \{1, 2, \dots
,n\}$.\end{ex}
\begin{table}
\begin{center}
    \begin{tabular}{| c || c | c | c |}\hline
    $T\sb 3$  &   1     &      2  &     3  \\ \hline \hline
            1 &   1     &      1  &     1   \\ \hline
            2 &   2     &      2  &     2   \\ \hline
            3 &   3     &      3  &     3   \\ \hline
    \end{tabular}
\caption{$T\sb 3$ Multiplication Table}\label{Ta:T3}
\end{center}
\end{table}
For $n=3$ its multiplication table is displayed in Table
\ref{Ta:T3}; the sequence of permutations is then $\big((), (),
()\big)$, where $()$ stands for the identity permutation. For
general $n$, the $i$-th term of the sequence of permutations is,
for any $i$ in $\{1, 2,\dots ,n\}$
\[
\mu\sb{i}=(\; )
\]

\begin{ex} $C\sb n$ is the {\bf cyclic rack} of order $n$. Its
underlying set is $\{1, 2, \dots ,n\}$ and $i\ast j := i+1$ mod
$n$, for any $i$, and $j\in \{1, 2, \dots ,n\}$.\end{ex}
\begin{table}
\begin{center}
    \begin{tabular}{| c || c | c | c |}\hline
    $C\sb 3$  &   1     &      2  &     3  \\ \hline \hline
            1 &   2     &      2  &     2   \\ \hline
            2 &   3     &      3  &     3   \\ \hline
            3 &   1     &      1  &     1   \\ \hline
    \end{tabular}
\caption{$C\sb 3$ Multiplication Table}\label{Ta:C3}
\end{center}
\end{table}
For $n=3$ its multiplication table is displayed in Table
\ref{Ta:C3}; the sequence of permutations is then $\big((1 \; 2 \;
3), (1 \; 2 \; 3), (1 \; 2 \; 3)\big)$. For general $n$, the
$i$-th term of the sequence of permutations is, for any $i$ in
$\{1, 2,\dots ,n\}$
\[
\mu\sb{i}=(\; 1 \; 2 \; 3 \; \cdots \; n\; )
\]

\section{Racks as Sequences of Permutations} \label{S:perm}

We begin with

\begin{def.} A {\bf rack} is a set, $X$, equipped with a binary operation, denoted $\ast$, such that:
\begin{itemize}
    \item for all $i$, and $j\in X$ there exists a unique $x\in X$ such that $x\ast j = i$ holds;
    \item for all $i$, $j$, and $k\in X$, $(i\ast j)\ast k=(i\ast k)\ast (j\ast k)$
\end{itemize}
\end{def.}

We refer to the first axiom as ``right-invertibility" and to the
second as ``self-distributivity".

\begin{def.} A {\bf quandle} is a rack, such that, for any $i\in X$, $i\ast i=i$. \end{def.}

This will be called the ``quandle condition".

\begin{def.} A {\bf prack} is a rack which is not a quandle.
\end{def.}

In this note, we will be dealing with finite racks or quandles.
Without loss of generality, the underlying set of such an
algebraic structure of order $n$, and any set of finite
cardinality $n$, will be taken to be the set of the first $n$
positive integers, $\{1, 2,\dots , n \}$, unless where explicitly
stated.

\begin{def.} A rack of finite order $n$ is called a {\bf $k$-prack} (for some $k$ in $\{0, 1, 2,\dots ,n\}$)
 if there are exactly $k$ distinct elements $i\sb {1}, i\sb {2},\dots
,i\sb {k}$ in $\{1, 2,\dots ,n\}$ such that $i\sb{j}\ast i\sb
{j}\neq i\sb {j}$.\end{def.}

Roughly speaking, a $k$-prack is a rack which violates $k$ times
the quandle condition. A $0$-prack is a quandle. In the
introduction, the trivial racks and the dihedral racks are
examples of quandles. The cyclic racks are examples of pracks; for
each $n$, the cyclic rack of order $n$ is an $n$-prack.

\begin{thm}\label{prop:condperm} Let $X$ be a finite set of cardinality $n$. Suppose a permutation $\mu\sb
i \in S\sb n$ (the group of permutations of the $n$ symbols $1,
2,\dots ,n$) is assigned to each $i$ in $X$. Then the expression
\[
\displaystyle i\ast j := \mu\sb{j}(i)
\]
for all $i$, $j \in \{1, 2,\dots ,n\}$, yields a rack structure in $X$ if and only if
\[
\displaystyle \mu\sb{\mu\sb{j}(i)}=\mu\sb j\mu\sb i\mu\sp{-1}\sb j
\]
for all $i$, $j \in \{1, 2,\dots ,n\}$. \end{thm}Proof: Suppose
$(X, \ast )$ is a rack. Then $(i\ast j)\ast k = (i\ast k)\ast
(j\ast k)$ is equivalent to $\mu\sb k \bigl(\mu\sb j
(i)\bigr)=\mu\sb{\mu\sb{k}(j)}\bigl( \mu\sb k (i)\bigr)$ which is
equivalent to $\mu\sb{\mu\sb{k}(j)}=\mu\sb k\mu\sb j\mu\sp{-1}\sb
k$. This proves the result in one direction. Now for the converse.
Self-distributivity is clear from the calculations above. Suppose
$i$, $j$, and $k \in \{1, 2,\dots ,n\}$. Then $k=i\ast j := \mu\sb
j(i)$ is equivalent to $\mu\sp{-1}\sb j(k)=i$. So, for any $k$, $j
\in \{1, 2,\dots ,n\}$  there is a unique solution to the equation
$k=x\ast j$. The result follows. \hfill $\blacksquare$

\begin{def.} Let $R=(X, \ast )$ be a finite rack of order $n$. We call {\bf sequence of permutations of $R$} the sequence
of permutations $(\mu\sb 1,\dots ,\mu\sb n)$ given by $\mu\sb j(i):=i\ast
j$, for all $i$, $j\in X$.
\end{def.}

In this way, Theorem \ref{prop:condperm} describes the conditions
the permutations of the sequence of permutations of a rack have to
satisfy. These permutations generate the {\bf operator group}
(\cite{dJoyce}, \cite{rFenn0}). We note

\begin{cor}Let $(\mu\sb 1,\dots ,\mu\sb n )$ denote the sequence of permutations of a rack $R$ of finite
order $n$. Then $R$ is a quandle if and only if $\mu\sb i(i)=i$,
for all $i$. \end{cor}Proof: Since $R=(X,\ast )$ is a rack then it
is a quandle if and only if $i=i\ast i = \mu\sb i(i)$, for all
$i$. \hfill $\blacksquare$

The right-invertibility axiom of a rack allows for the definition
of a second binary operation on the rack.

\begin{def.} Let $R=(X, \ast )$ be a rack. For any $i$, $j\in X$, let $i\:\bar{\ast}\:j$ denote the
 unique solution to the
equation $x\ast j =i$. $\:\bar{\ast}\:$ is said to be the dual
operation to $\ast $. $X$ equipped with this binary operation is
called the {\bf dual} of $R$.
\end{def.}

In the proof of Theorem \ref{prop:condperm} we came up with a
formula for this unique solution provided we had the original rack
represented by its sequence of permutations. This will help to
prove that

\begin{cor}The dual of a rack (resp., quandle) is a rack (resp., quandle). \end{cor}Proof: Let $R=(X,\ast )$
 be a finite rack with sequence of permutations $(\mu\sb 1,\dots ,\mu\sb n)$. As seen above, the dual operation is
given by $i\:\bar{\ast}\: j = \mu\sp{-1}\sb j(i)$. So if the dual
of $R$ is also a rack, its sequence of permutations has to be
$(\mu\sp{-1}\sb 1,\dots ,\mu\sp{-1}\sb n)$. Since $R$ is a rack
then, for all $i$, $j$, $\mu\sb{\mu\sb{i}(j)}=\mu\sb i\mu\sb
j\mu\sb i\sp{-1}$ which, by taking inverses, is equivalent to
$\mu\sp{-1}\sb{\mu\sb{i}(j)}=\mu\sb i\mu\sp{-1}\sb j\mu\sp{-1}\sb
i$ which can be rewritten as $\mu\sp{-1}\sb
i\mu\sp{-1}\sb{\mu\sb{i}(j)}\mu\sb i=\mu\sp{-1}\sb j$. Since, for
any $i$, $\mu\sb i$ is a bijection, then we can put $j'=\mu\sb
i(j)$ so that $j=\mu\sp{-1}\sb i(j')$ and so, for all $j'$ and
$i$, $\mu\sp{-1}\sb i\mu\sp{-1}\sb{j'}\mu\sb i=\mu\sp{-1}\sb
{\mu\sp{-1}\sb{i}(j')}$. Finally, if $R$ is a quandle then, for
any $i$, $\mu\sb{i}(i)=i$ which is equivalent to saying that
$\mu\sb{i}\sp{-1}(i)=i$, for all $i$. The result follows.
$\blacksquare$

\begin{cor}The $\ast$ and $\:\bar{\ast}\:$ operations are distributive with respect to each other
i.e., for all $i$, $j$ and $k$,
\[
(i\ast j)\:\bar{\ast}\: k = (i\:\bar{\ast}\: k)\ast
(j\:\bar{\ast}\: k) \qquad  and \qquad (i\:\bar{\ast}\: j)\ast k =
(i\ast k)\:\bar{\ast}\: (j\ast k)
\]\end{cor}
Proof: We will just prove that $(i\ast j)\:\bar{\ast}\: k =
(i\:\bar{\ast}\: k)\ast (j\:\bar{\ast}\: k)$. This is equivalent
to
$\mu\sp{-1}\sb{k}\mu\sb{j}(i)=\mu\sb{\mu\sp{-1}(j)}\mu\sp{-1}\mu\sb{k}(i)$
which is equivalent to
$\mu\sp{-1}\sb{k}\mu\sp{-1}\sb{j}\mu\sb{k}=\mu\sp{-1}\sb{\mu\sp{-1}\sb{k}(j)}$
which is the statement that the $\:\bar{\ast}\:$ is
self-distributive. $\blacksquare$

\begin{prop}\label{prop:musp{n}1} Let R be a finite rack of order $n$. There exists a positive
integer $l$ such that, for any $j$ in $\{1, 2,\dots , n \}$,
$\mu\sb j\sp l = ()$, where $()$ stands for the identity
permutation.
\end{prop}Proof: Write the $\mu\sb{i}$'s as products of disjoint cycles and
let $\Sigma$ stand for the set of all these cycles, except
$1$-cycles. Then $l$ can be taken to be the least common
multiplier of the lengths of the cycles in $\Sigma$.
$\blacksquare$

For the purpose of the discussion in this section we will fix two
racks, $R$ and $R'$, of finite order $n$ but otherwise arbitrary.
The operations will be denoted $\ast$ and $\ast'$, respectively.
The sequences of their permutations will be denoted $(\mu\sb
1,\dots ,\mu\sb n)\subset S\sb n$, and $(\mu'\sb 1,\dots ,\mu'\sb
n)\subset S\sb n$, respectively.

\begin{def.} A {\bf rack homomorphism} between two racks, $R$ and $R'$, is a map, $\alpha $,
 between the underlying sets, such that $\alpha (i\ast j)=\alpha(i)\ast '\alpha(j)$,
for all $i$, and $j$. If $\alpha$ is a bijection then it is a {\bf
rack isomorphism}.
\end{def.}

Clearly, if two racks are isomorphic they have to have the same
order. Moreover,

\begin{prop}\label{prop:qnr}If a $k$-prack, $R$, is isomorphic to a $k'$-prack, $R'$, then $k=k'$.\end{prop}
Proof: Assume to the contrary and let $k < k'$. Let $\alpha $ be
the isomorphism from the $k'$-prack to the $k$-prack. Let
$i\sb{1},\dots ,i\sb{k'}$ be all elements of the $k'$-prack which
violate the ``quandle condition''. Then $\alpha (i\sb{j})\ast
'\alpha (i\sb{j}) \neq \alpha (i\sb{j})$ hence the image cannot be
a $k$-prack.$\blacksquare$

Furthermore,

\begin{thm}\label{prop:rackiso}$\alpha $ is a rack isomorphism of $R$ with $R'$ if and only if:
\[
\displaystyle \mu'\sb{\alpha(j)}=\alpha \mu\sb{j}\alpha\sp{-1}
\]
for each $j \in \{1, 2,\dots ,n\}$.
\end{thm}Proof: For any $i$, and $j\in X$, $\alpha (i\ast j)=\alpha(i)\ast ' \alpha(j)$, which, in terms
of the sequences of permutations is $\alpha
(\mu\sb{j}(i))=\mu'\sb{\alpha(j)}\bigl( \alpha(i)\bigr)$, and
since $\alpha$ is a bijection:
\[
\displaystyle \mu'\sb{\alpha(j)}=\alpha \mu\sb{j}\alpha\sp{-1}
\]
and conversely. $\blacksquare$

In this way a rack isomorphism $\alpha$, being a permutation from
$S\sb n$, gives rise to an action on $S\sb n$ by conjugation
\[
\begin{array}{clcr}
\alpha   \longmapsto \hat{\alpha} : & S\sb n \longrightarrow S\sb n\\
        & \mu \longmapsto \alpha\mu \alpha\sp{-1}
\end{array}
\]
which induces a bijection between the sequences of permutations of
the two racks, according to:
\[
\displaystyle \mu'\sb{\alpha(j)}=\alpha \mu\sb{j}\alpha\sp{-1}
\]

\section{Invariants of Racks - Basics} \label{S:inv1}

Theorem \ref{prop:rackiso} gives necessary conditions on two
sequences of permutations in order for the corresponding racks to
be isomorphic. Before looking at some examples, here is a

\begin{def.}Let $\nu \in S\sb n$ be a permutation whose decomposition into disjoint cycles
is (modulo ordering of the cycles),
\[
\nu = (i\sp 1\sb 1,\dots ,i\sp 1\sb{n\sb {1}})(i\sp 2\sb 1,\dots ,i\sp 2\sb {n\sb {2}})\dots (i\sp k\sb 1,\dots ,i\sp k\sb {n\sb {k}})
\]
We call the list of the lengths of the cycles, $\{n\sb 1,\dots ,n\sb k\}$, the {\bf pattern} of $\nu$.
\end{def.}
with the following understanding,
\begin{def.} A list is a set with (possible) repeats.\end{def.}

The notation for list will be the same as for sets i.e., the
elements will be involved by curled braces. It will be clear
either from context or by explicit mention which object we mean.

Consider the following five racks of order four given by their sequences of permutations where the
permutations are broken down into products of disjoint cycles.
\[
R: \quad \bigl( (), (), (1 2)(3)(4), (1 2)(3)(4)\bigr)
\]
\[
S: \quad \bigl( (), (), (1)(2)(3 4), (1)(2)(3 4)\bigr)
\]
\[
T: \quad \bigl( (1 2 3)(4), (1 2 3)(4), (1 2 3)(4), (1 2
3)(4)\bigr)
\]
\[
U: \quad \bigl( (1 2 3)(4), (1 2 3)(4), (1 2 3)(4), (1 3
2)(4)\bigr)
\]
\[
V: \quad \bigl( (1)(2 3 4), (2)(1 4 3), (3)(1 2 4), (4)(1 3
2)\bigr)
\]
Clearly, none of the first two is isomorphic to any of the last
three racks since conjugation preserves the pattern of the
permutations. In this way, we define a first invariant of rack
isomorphism.

\begin{def.}Let $R$ be a rack with sequence of permutations $(\mu\sb 1,\dots ,\mu\sb n)\subset S\sb n$.
Let  $\{ k\sp i\sb 1,\dots ,k\sp i\sb{m\sb i}\}$ be the pattern of $\mu\sb i$ for
each $i\in \{1,\dots ,n\}$. We call the sequence of patterns
$\bigl( \{ k\sp 1\sb 1,\dots ,k\sp 1\sb{m\sb 1}\} ,\dots ,\{ k\sp
n\sb 1,\dots ,k\sp n\sb{m\sb n}\} \bigr)$ the {\bf profile} of
$R$.
\end{def.}

\begin{prop}\label{prop:profile} The profile is an invariant under rack isomorphism. \end{prop}
Proof: From Theorem \ref{prop:rackiso} and remark right after it,
any rack isomorphism establishes a one-to-one correspondence
between the sequences of permutations of the racks under study via
conjugation by a fixed permutation. Since conjugation preserves
the pattern the result follows.  $\blacksquare$

We can now subdivide these five racks, $R$, $S$, $T$, $U$, and
$V$, into two classes according to profile (see Table
\ref{Ta:RP}).
\begin{table}
\begin{center}
\renewcommand{\arraystretch}{1.25}
    \begin{tabular}{| c | c | c |}\hline
Racks     &             Profile & Sequences of Permutations
\\ \hline R         & $\bigl($\{1, 1, 1, 1\}, \{1, 1, 1, 1\}, \{1,
1, 2\}, \{1, 1, 2\}$\bigr)$ & ((), (), (1 2)(3)(4), (1 2)(3)(4))
\\\cline{1-1}  \cline{3-3} S & &     ((), (), (1)(2)(3 4), (1)(2)(3 4))\\
\hline
T         & &     $\bigl($(1 2 3)(4), (1 2 3)(4), (1 2 3)(4), (1 2 3)(4)$\bigr)$            \\
\cline{1-1} \cline{3-3} U         &        $\bigl($\{1, 3\}, \{1,
3\}, \{1, 3\}, \{1, 3\}$\bigr)$                         &     ((1
2
3)(4), (1 2 3)(4), (1 2 3)(4), (1 3 2)(4))            \\
\cline{1-1}\cline{3-3} V & & ((1)(2 3 4), (1 4 3)(2), (1 2 4)(3),
(1 3 2)(4))
\\  \hline
    \end{tabular}
    \caption{Racks and Their Profiles}\label{Ta:RP}
\end{center}
\end{table}
Note that, for the first profile, $\bigl(\{1, 1, 1, 1\}, \{1, 1,
1, 1\}, \{1, 1, 2\}, \{1, 1, 2\}\bigr)$, there are two distinct
patterns, $\{1, 1, 1, 1\}$ and $\{1, 1, 2\}$, each of them
realized in just one way with multiplicity two, for each rack (see
Table \ref{Ta:PR}).
\begin{table}
\begin{center}
\renewcommand{\arraystretch}{1.25}
    \begin{tabular}{| c | c | c | c |}\hline
Racks     &       Profile
&          Distinct Patterns          &      ...realized as...  \\
\hline R     &
&      \{1, 1, 1, 1\}                   &     (), ()            \\
\cline{3-4}
          &       $\bigl($\{1, 1, 1, 1\}, \{1, 1, 1, 1\}, \{1, 1, 2\}, \{1, 1, 2\}$\bigr)$                  &       \{1, 1, 2\}                     &     (1 2)(3)(4), (1 2)(3)(4)        \\ \cline{1-1} \cline{3-4}
S     &
&      \{1, 1, 1, 1\}                   &     (), ()            \\
\cline{3-4}
          &                                                                             &        \{1, 1, 2\}                    &     (1)(2)(3 4), (1)(2)(3 4)        \\ \hline
    \end{tabular}
    \caption{How the profile is realized for $R$ and $S$}\label{Ta:PR}
\end{center}
\end{table}
As for the second profile, $\bigl(\{1, 3\}, \{1, 3\}, \{1, 3\},
\{1, 3\}\bigr)$, there is only one distinct pattern which is
realized in different ways for each of the different racks, $T$,
$U$, and $V$ (see Table \ref{Ta:PPR}).
\begin{table}
\begin{center}
\renewcommand{\arraystretch}{1.25}
    \begin{tabular}{| c | c | c | c |}\hline
Racks     &       Profile                                       &
Distinct Patterns          &      ...realized as...             \\
\hline T     & &      \{1, 3\}                         &     (1 2
3)(4), (1 2 3)(4), (1 2 3)(4), (1 2 3)(4)            \\
\cline{1-1} \cline{3-4} U         & &       \{1, 3\} & (1 2 3)(4),
(1 2 3)(4), (1 2 3)(4)        \\ \cline{4-4}
          &   $\bigl($\{1, 3\}, \{1, 3\}, \{1, 3\}, \{1, 3\}$\bigr)$                &                                     &     (1 3 2)(4)                    \\  \cline{1-1} \cline{3-4}
V     &                                                     &
\{1, 3\}                         &     (1)(2 3 4)            \\
\cline{4-4}
          &                                                     &                                     &     (1 4 3)(2)        \\   \cline{4-4}
      &                                                     &                                     &     (1 2 4)(3)            \\ \cline{4-4}
          &                                                     &                                     &     (1 3 2)(4)        \\   \hline
    \end{tabular}
    \caption{How the profile is realized for $T$, $U$ and $V$}\label{Ta:PPR}
\end{center}
\end{table}
Now $T$ cannot be isomorphic either to $U$ or to $V$ because the isomorphism would send the four equal
permutations into another four equal permutations and neither $U$ nor $V$ have such a
sequence of permutations. Analogously $U$ cannot be isomorphic to the other two and $V$ cannot be
isomorphic to the other two.

In this way we give the following

\begin{def.}Suppose $R$ is a rack of finite order $n$ with sequence of permutations
$(\mu\sb 1,\dots ,\mu\sb n)\subset S\sb n$ and profile $\bigl( \{n\sp 1\sb
1,\dots ,n\sp 1\sb{k\sb 1}\} ,\dots ,\{ n\sp n\sb 1,\dots ,n\sp
n\sb{n\sb{k\sb n}}\} \bigr)$.

Assume further that the sequence of {\bf distinct} patterns is
$\bigl( \{n\sp 1\sb 1,\dots ,n\sp 1\sb{n\sb 1}\},\dots ,\{n\sp
m\sb 1,\dots ,n\sp m\sb{n\sb m}\}\bigr)$ with $m\leq n$, with
$r\sp i$ repeats for the $i$-th distinct pattern,
 corresponding to $r\sp i\sb 1$ identical permutations, $\rho\sp i\sb 1$, ..., $r\sp i\sb{j\sb i}$ identical permutations,
 $\rho\sp i\sb{j\sb i}$, such that
 $r\sp i\sb 1+\dots +r\sp i\sb{j\sb i}=r\sp i$ and $\rho\sp i\sb k\neq \rho\sp i\sb l$ for $k\neq l$.

For each of the distinct patterns, say $\{n\sp i\sb 1,\dots ,n\sp
i\sb{n\sb i}\}$ consider the ordered pair $\bigl( \{n\sp i\sb
1,\dots ,n\sp i\sb{n\sb i}\}, \{r\sp i\sb 1,\dots ,r\sp i\sb{n\sb
i}\} \bigr)$, whose first term is this pattern and whose second
term is the list of multiplicities corresponding to this pattern,
$\{r\sp i\sb 1,\dots ,r\sp i\sb{n\sb i}\}$.

We call the sequence of these ordered pairs, one ordered pair per
distinct pattern,
\[
\Bigl( \bigl( \{n\sb 1\sp 1,\dots ,n\sp 1\sb{n\sb 1}\},\{r\sp 1\sb
1,\dots ,r\sp 1\sb{j\sb 1}\}\bigr),\dots ,\bigl( \{n\sb 1\sp
m,\dots ,n\sp m\sb{n\sb m}\},\{r\sp m\sb 1,\dots ,r\sp m\sb{j\sb
m}\}\bigr)\Bigr)
\]
the {\bf detail} of $R$.\end{def.}

\begin{def.}Keeping the notation of the previous definition, we
call the sequence
\[
\Bigl(\bigl(\{n\sb 1\sp 1,\dots ,n\sp 1\sb{n\sb 1}\},\{(r\sp 1\sb
1, \rho\sp 1\sb 1),\dots ,(r\sp 1\sb{j\sb 1}, \rho\sp 1\sb{j\sb
1}\}\bigr),\dots ,\bigl(\{n\sb 1\sp m,\dots ,n\sp m\sb{n\sb
m}\},\{(r\sp m\sb 1, \rho\sp m\sb 1),\dots ,(r\sp m\sb{j\sb m},
\rho\sp m\sb{j\sb m})\}\bigr)\Bigr)
\]
the {\bf absolute detail} of $R$.\end{def.}

\begin{prop}\label{prop:detail} The detail is an invariant under rack isomorphism. \end{prop}
Proof: Omitted since it is similar to the proof of the previous
proposition.  $\blacksquare$

We can now break down the profile $\bigl(\{1, 3\}, \{1, 3\}, \{1,
3\}, \{1, 3\}\bigr)$ into three distinct details (see Table
\ref{Ta:DP}).
\begin{table}
\begin{center}
\renewcommand{\arraystretch}{2.0}
    \begin{tabular}{| c | c | c |}\hline
Rack       &     Detail                           & Absolute
Detail     \\ \hline T         &     $\Bigl( \bigl($\{1 ,3\},
\{4\}$\bigr) \Bigr)$ & $\Bigl( \bigl($\{1, 3\}, \{(4, (1 2 3)(4))\}$\bigr) \Bigr)$  \\
\cline{1-1} \cline{2-3} U & $\bigl( \bigl($\{1 ,3\}, \{1,
3\}$\bigr) \bigr)$ & $\Bigl( \bigl($\{1, 3\}, \{(1, (1 3 2)(4)),
(3, (1 2 3)(4))\}$\bigr) \Bigr)$
\\ \cline{1-1} \cline{2-3} V         &  $\Bigl( \bigl($\{1 ,3\}, \{1, 1, 1, 1\}$\bigr)
\Bigr)$ &         $\Bigl( \bigl($\{1, 3\}, \{(1, (1)(2 3 4)), (1,
(1 4 3)(2)), (1, (1 2 4)(3)), (1, (1 3 2)(4))\}$\bigr) \Bigr)$ \\
\hline
    \end{tabular}
    \caption{Profile $\{(1, 3), (1, 3), (1, 3), (1, 3)\}$ has different details}\label{Ta:DP}
\end{center}
\end{table}

\begin{cor} Isomorphic racks have the same detail. Racks with the same detail have the same profile.\end{cor}
Proof: Clear from the discussion above. $\blacksquare$

The converses are not true. Specifically, we saw that racks $T$,
$U$, and $V$ have the same profile but different details, and $R$
and $S$ although they have the same detail, they are not
isomorphic, since $R$ is a quandle and $S$ is a $2$-prack.

\section{Invariants of Racks - More Examples} \label{S:inv2}

From the considerations of the previous Section we can extract a
sequence of instructions for a partial classification of racks for
each order $n$. Upon selection of the order $n$, all racks of that
order are generated, each one being specified by its sequence of
permutations. Then they are sorted according to their profiles.
Then within each profile, they are sorted according to their
detail. Then, within each detail, they are sorted according to
being a $k$-prack for each $k$ in $\{0, 1, 2,\dots ,n\}$. Of
course, we still end up with a number of different but isomorphic
racks. That is why we now introduce some more invariants of racks
by discussing examples.

The Laurent polynomials have a natural structure of quandle given
by setting $a\ast b:=Ta+(1-T)b$ for any $a$, $b\in \mathbb{Z}[T,
T\sp{-1}]$. This is an infinite quandle. In order to obtain finite
quandles, we take the coefficients of the polynomials from
$\mathbb{Z}\sb{n}$ and take the polynomials from the quotient of
$\mathbb{Z}\sb{n}[T, T\sp{-1}]$ by some principal ideal, being
careful in the choice of $n$ and ideal so that the
right-invertibility axiom holds. These are the so-called finite
Alexander quandles, denoted $\mathbb{Z}\sb{n}[T,
T\sp{-1}]/(h(T))$, where $(h(T))$ is the principal ideal generated
by polynomial $h$. For each $n$ greater than $2$, the dihedral
quandle $R\sb n$ is isomorphic with the finite Alexander quandle
$\mathbb{Z}\sb{n}[T, T\sp{-1}]/(T+1)$ (\cite{dJoyce},
\cite{sNelson}).

 It is known that, for any
prime order $p \geq 3$, there are $p-2$ non-trivial quandles
pairwise non-isomorphic, which are isomorphic with the linear
Alexander quandles of order $p$ i.e., $\mathbb{Z}\sb{p}[T,
T\sp{-1}]/(T-m)$ where $m$ is from $\{2, 3, \dots ,p-1 \}$ and
$a\ast b = Ta+(1-T)b$ in the indicated quotient (\cite{sNelson}).
We will now concentrate on how to tell these racks apart for $p=5$
using the techniques developed in the previous Section. There are
three non-trivial linear Alexander quandles of order $5$,
$\mathbb{Z}\sb{5}[T, T\sp{-1}]/(T-2)$, $\mathbb{Z}\sb{5}[T,
T\sp{-1}]/(T-3)$, and $\mathbb{Z}\sb{5}[T, T\sp{-1}]/(T-4)$. The
sequences of permutations are, respectively:
\[
\mu\sb{1}=(1)(2354), \quad \mu\sb{2}=(2)(1534), \quad
\mu\sb{3}=(3)(1452), \quad \mu\sb{4}=(4)(1325), \quad
\mu\sb{5}=(5)(1243)
\]
(see multiplication table in Table \ref{Ta:m=2} where, for
instance $2\ast 1 = 3$)
\begin{table}[h!]
\begin{center}
    \begin{tabular}{| c || c | c | c | c | c |}\hline
              &   1     &      2  &     3   &    4   &   5   \\ \hline \hline
            1 &   1     &      5  &     4   &    3   &   2      \\ \hline
            2 &   3     &      2  &     1   &    5   &   4      \\ \hline
            3 &   5     &      4  &     3   &    2   &   1      \\ \hline
            4 &   2     &      1  &     5   &    4   &   3      \\ \hline
            5 &   4     &      3  &     2   &    1   &   5      \\ \hline
    \end{tabular}
\caption{$\mathbb{Z}\sb{5}[T, T\sp{-1}]/(T-2)$ Multiplication
Table}\label{Ta:m=2}
\end{center}
\end{table},
\[
\mu\sb{1}=(1)(2453), \quad \mu\sb{2}=(2)(1435), \quad \mu\sb{3}=(3)(1254), \quad \mu\sb{4}=(4)(1523), \quad \mu\sb{5}=(5)(1342)
\]
(see multiplication table in Table \ref{Ta:m=3})
\begin{table}[h!]
\begin{center}
    \begin{tabular}{| c || c | c | c | c | c |}\hline
              &   1     &      2  &     3   &    4   &   5   \\ \hline \hline
            1 &   1     &      4  &     2   &    5   &   3      \\ \hline
            2 &   4     &      2  &     5   &    3   &   1      \\ \hline
            3 &   2     &      5  &     3   &    1   &   4      \\ \hline
            4 &   5     &      3  &     1   &    4   &   2      \\ \hline
            5 &   3     &      1  &     4   &    2   &   5      \\ \hline
    \end{tabular}
\caption{$\mathbb{Z}\sb{5}[T, T\sp{-1}]/(T-3)$ Multiplication
Table}\label{Ta:m=3}
\end{center}
\end{table}, and,
\[
\mu\sb{1}=(1)(25)(34), \quad \mu\sb{2}=(2)(13)(45), \quad \mu\sb{3}=(3)(15)(24), \quad \mu\sb{4}=(4)(12)(35), \quad \mu\sb{5}=(5)(14)(23)
\]
(see multiplication table in Table \ref{Ta:m=4}).
\begin{table}[h!]
\begin{center}
    \begin{tabular}{| c || c | c | c | c | c |}\hline
              &   1     &      2  &     3   &    4   &   5   \\ \hline \hline
            1 &   1     &      3  &     5   &    2   &   4      \\ \hline
            2 &   5     &      2  &     4   &    1   &   3      \\ \hline
            3 &   4     &      1  &     3   &    5   &   2      \\ \hline
            4 &   3     &      5  &     2   &    4   &   1      \\ \hline
            5 &   2     &      4  &     1   &    3   &   5      \\ \hline
    \end{tabular}
\caption{$\mathbb{Z}\sb{5}[T, T\sp{-1}]/(T-4)$ Multiplication
Table}\label{Ta:m=4}
\end{center}
\end{table}

Now $\mathbb{Z}\sb{5}[T, T\sp{-1}]/(T-4)$ is not isomorphic with
$\mathbb{Z}\sb{5}[T, T\sp{-1}]/(T-3)$ nor with
$\mathbb{Z}\sb{5}[T, T\sp{-1}]/(T-2)$ by inspection of profiles.
As for $\mathbb{Z}\sb{5}[T, T\sp{-1}]/(T-3)$ and
$\mathbb{Z}\sb{5}[T, T\sp{-1}]/(T-2)$, they have the same detail.
Inspection of their multiplication tables shows that that of
$\mathbb{Z}\sb{5}[T, T\sp{-1}]/(T-3)$ is symmetric with respect to
the diagonal whereas that of $\mathbb{Z}\sb{5}[T, T\sp{-1}]/(T-2)$
is not. Clearly, an isomorphism would preserve this symmetry hence
$\mathbb{Z}\sb{5}[T, T\sp{-1}]/(T-2)$ and $\mathbb{Z}\sb{5}[T,
T\sp{-1}]/(T-3)$ are not isomorphic. We further remark that these
two quandles form an example of a rack which is not isomorphic to
its dual. In fact, the sequence of permutations, $(\nu\sb{1},
\nu\sb{2}, \nu\sb{3}, \nu\sb{4}, \nu\sb{5})$, of the dual of
$\mathbb{Z}\sb{5}[T, T\sp{-1}]/(T-3)$ is, omitting $1$-cycles,
\begin{align}\notag
\nu\sb{1} & =(2453)\sp{-1} = (3542) , \quad \nu\sb{2} =(1435)
\sp{-1}  = (5341) , \quad \nu\sb{3}= (1254)\sp{-1} = (4521) ,\\
\notag
 \quad \nu\sb{4} & =(1523) \sp{-1}
= (3251) , \quad \nu\sb{5} = (1342)\sp{-1} = (2431)
\end{align}
which is the sequence of permutations of $\mathbb{Z}\sb{5}[T,
T\sp{-1}]/(T-2)$. Since these two racks are not isomorphic the
remark follows.

We will now develop the notion of symmetry. Let $R=(X, \ast)$ and
$(\mu\sb{1},\dots ,\mu\sb{n})$ have the usual meanings.

\begin{def.}Given a rack, $R$, we call {\bf opposite} of $R$, notation $R\sp{o}=(X, \ast \sp{o})$,
the algebraic structure whose underlying set is again $X=\{1, 2,\dots ,n\}$ and whose binary operation
is given by:
\[
i\ast\sp{o} j := j\ast i = \mu\sb{i}(j)
\]
for any $i$, and $j$ in $\{1, 2,\dots ,n\}$.\end{def.}

\begin{def.}Keeping the notation above, a rack is called {\bf symmetric} if for any $i$, and
$j$ in $\{1, 2,\dots ,n\}$, $i\ast\sp{o} j = i\ast j$. \end{def.}

\begin{prop} Symmetry is an invariant under rack isomorphism.\end{prop}Proof: Let $\alpha$
 denote the isomorphism. Then
 $\alpha (i){\ast'}\sp{o} \alpha (j)=\alpha (j)\ast' \alpha (i)=\alpha (j\ast i)=\alpha (i\ast\sp{o}j)=\alpha (i\ast j)=\alpha (i)\ast'\alpha (j)$. $\blacksquare$

When is the opposite of a rack (resp., quandle) a rack (resp.
quandle)? Clearly,

\begin{prop}\label{prop: tq} If $R=(X, \ast )$ is a quandle, then all elements in $R\sp{o}$ satisfy the
quandle condition i.e., $i\ast\sp{o}i=i$ for all
$i$.\end{prop}Proof: Omitted. $\blacksquare$

What about right-invertibility?

\begin{prop}\label{prop: 1-c} Let $R=(X, \ast )$ be a finite rack.  $\ast\sp{o}$ is right-invertible if and only if, for any
$i$, and $j$ in $\{1, 2,\dots ,n\}$, there exists a unique $x$ in
$\{1, 2,\dots ,n\}$ such that $i\ast x=j$.\end{prop}Proof: If
$R\sp{o}$ is right-invertible then, for any $i$, and $j$ in $\{1,
2,\dots ,n\}$, there exists a unique $x$ in $\{1, 2,\dots ,n\}$
such that $j=x\ast\sp{o}i=i\ast x$ and vice-versa. $\blacksquare$

The racks that satisfy the condition on Proposition \ref{prop:
1-c} are both left- and right-invertible. An algebraic structure
whose binary operation is both left- and right-invertible is
called a quasi-group (\cite{SR}). Clearly, $\mathbb{Z}\sb{5}[T,
T\sp{-1}]/(T-3)$ is a quasi-group since being symmetric, its
opposite is again a quandle. Furthermore, dihedral quandles of any
odd order, say $2m+1$, are quasi-groups since for any $i$, and $j$
in $\{1, 2,\dots ,2m+1 \}$ the equation $i\ast x = j$ is
equivalent to $2x = i+j$, which, upon multiplication on both sides
by $m+1$, yields the unique $x=(m+1)(i+j)$. Straightforward
calculations show also that for each prime order $p$ the
non-trivial (linear) Alexander quandles are quasi-groups. We will
now give a partial answer to the question of when the opposite of
a quandle is again a quandle. We start out by noting that

\begin{prop}\label{prop: aqsd} If $(X, \ast )$ is an Alexander quandle, then its opposite is self-distributive. \end{prop}
Proof: What we have to prove is that $(i\ast\sp{o}j)\ast\sp{o}k =
(i\ast\sp{o}k)\ast\sp{o}(j\ast\sp{o}k)$, that is, $k\ast(j\ast i)
= (k\ast j)\ast (k\ast i)$, that is, $Tk+(1-T)\bigl(
Tj+(1-T)i\bigr) =T\bigl( Tk+(1-T)j\bigr)+(1-T)\bigl(
Tk+(1-T)i\bigr) $ and this is now a straightforward calculation.
$\blacksquare$

In this way, Propositions \ref{prop: tq} and \ref{prop: aqsd} show
that the opposite of any Alexander quandle satisfies the quandle
axiom and the self-distributivity axiom. From the remarks made
right before Proposition \ref{prop: aqsd} we can say that for any
odd order, say $2m+1$, the opposite of the dihedral quandle
$R\sb{2m+1}$ is again a quandle; and that the opposite of any non
trivial linear Alexander quandle of any prime order is again a
quandle.

We will now give an example of a quandle whose opposite is not
self-distributive. In Table \ref{Ta:octa} we see the
multiplication table of a quandle such that its opposite is not
right-invertible. This is the octahedron quandle (\cite{AG}).
\begin{table}[h!]
\begin{center}
    \begin{tabular}{| c || c | c | c | c | c | c |}\hline
              &   1     &      2  &     3   &    4   &   5   &   6   \\ \hline \hline
            1 &   1     &      5  &     2   &    3   &   4   &   1   \\ \hline
            2 &   3     &      2  &     6   &    2   &   1   &   5   \\ \hline
            3 &   4     &      1  &     3   &    6   &   3   &   2   \\ \hline
            4 &   5     &      4  &     1   &    4   &   6   &   3   \\ \hline
            5 &   2     &      6  &     5   &    1   &   5   &   4   \\ \hline
            6 &   6     &      3  &     4   &    5   &   2   &   6   \\ \hline
    \end{tabular}
\caption{Octahedron Quandle Multiplication Table}\label{Ta:octa}
\end{center}
\end{table}
Omitting $1$-cycles, its sequence of permutations is
\[
\mu\sb{1}=(2345), \quad \mu\sb{2}=(1563), \quad \mu\sb{3}=(1264), \quad \mu\sb{4}=(1365), \quad \mu\sb{5}=(1462), \quad \mu\sb{6}=(2543)
\]

We remark that, for the octahedron quandle, we have $1\ast (2\ast
3)=1\ast 6=1$ whereas $(1\ast 2)\ast (1\ast 3)=5\ast 2=6$, thus
the opposite of the octahedron quandle is not self-distributive.
Hence, the octahedron quandle is not an Alexander quandle.

In Proposition \ref{prop: 1-c} we dealt with this property of
having just one element connecting any two elements of the rack.
This prompted us into the following. Let $R=(X, \ast )$ be a
finite rack. For any $i$, and $j$ in $\{1, 2,\dots ,n\}$ set:
\[
N\sb{i, j}:=\{m\in \mathbb{Z}\sb{0}\sp{+} \; | \; \exists \;
k\sb{1},\dots ,k\sb{m}\in \{1, 2,\dots ,n\}, \;
\epsilon\sb{1},\dots ,\epsilon\sb{m} \in \{ \pm 1 \}  \; : \;
j=(\cdots ((i\ast \sp{\epsilon
\sb{1}}k\sb{1})\ast\sp{\epsilon\sb{2}}k\sb{2})\cdots
)\ast\sp{\epsilon\sb{m}}k\sb{m} \}
\]
(where $\ast \sp{+1}:=\ast $ and $\ast \sp{-1}:=\:\bar{\ast}\:$).
If $N\sb{i, j}\neq \emptyset$ set $m\sb{i, j}:=\min N\sb{i, j}$, otherwise $m\sb{i, j}:=0$.

\begin{def.}Keeping the notation above, we say a rack is {\bf $c$-connected} if, for all $i$, and $j$ in
$\{1, 2,\dots ,n\}$, $m(i, j)\neq \emptyset$ and
\[
c=\max \{ m(i, j) \: |\: i, j \in \{1, 2,\dots ,n\} \}
\]
Otherwise, if there are $i$, and $j$ in $\{1, 2,\dots ,n\}$, such
that $m(i, j)= \emptyset$, then we say the rack is {\bf not
connected}.
\end{def.}

In this way, the octahedron quandle is $2$-connected.

\begin{prop} Keeping the notation above, if a rack of finite order $n$ is $c$-connected for $c>0$, then $c\leq n$.
\end{prop}Proof: Omitted. $\blacksquare$

We will now prove that $c$-connectedness is an invariant of racks
under isomorphism after introducing the relevant object for this
proof.

\begin{def.} Let  $R=(X, \ast )$ be a finite rack. For each $i$ in $\{1, 2,\dots ,n\}$ we call the set:
\[
{\cal O}\sb{i}:=\{ (\cdots ((i\ast\sp{\epsilon\sb{1}}k\sb{1})\ast\sp{\epsilon\sb{2}}k\sb{2})\cdots )\ast\sp{\epsilon\sb{m}}k\sb{m} \; | \; m\in \mathbb{Z}\sb{0}\sp{+}, \:
k\sb{1},\dots ,k\sb{m} \in \{1,\dots ,n \}, \: \epsilon\sb{1},\dots
\epsilon\sb{m} \in \{ \pm 1 \} \}
\]
the {\bf orbit of $i$}.
\end{def.}

\begin{prop} Let $\alpha$ be a rack isomorphism between racks $R$ and $R'$. For any $i$, let ${\cal O}\sb{i}$ and
${\cal O'}\sb{i}$ denote the orbit of $i$ in $R$ and $R'$,
respectively. Then $\alpha$ restricted to ${\cal O}\sb{i}$ is a
bijection of ${\cal O}\sb{i}$ with ${\cal O'}\sb{\alpha (i)}$.
Furthermore, $\alpha$ induces a bijection from the set of orbits
of $R$ to the set of orbits of $R'$. \end{prop}Proof: We just
remark that
\begin{align}\notag
& \alpha \Bigl( {\cal O}\sb{i}\Bigr)  =\alpha \Big( \{ (\cdots
((i\ast\sp{\epsilon\sb{1}}k\sb{1})\ast\sp{\epsilon\sb{2}}k\sb{2})\cdots
)\ast\sp{\epsilon\sb{m}}k\sb{m} \; | \; m\in
\mathbb{Z}\sb{0}\sp{+}, \: k\sb{1},\dots ,k\sb{m} \in \{1,\dots ,n
\}, \: \epsilon\sb{1},\dots \epsilon\sb{m} \in \{ \pm 1 \}
\}\Bigr)= \\ \notag
 & = \{ (\cdots ((\alpha (i){\ast'}\sp{\epsilon\sb{1}}\alpha
(k\sb{1})){\ast'}\sp{\epsilon\sb{2}}\alpha (k\sb{2}))\cdots
){\ast'}\sp{\epsilon\sb{m}}\alpha (k\sb{m}) \; | \; m\in
\mathbb{Z}\sb{0}\sp{+}, \: k\sb{1},\dots ,k\sb{m} \in \{1,\dots ,n
\}, \: \epsilon\sb{1},\dots \epsilon\sb{m} \in \{ \pm 1 \} \}
\end{align}
$\blacksquare$

\begin{prop} $c$-connectedness is an invariant under rack isomorphism. \end{prop}Proof: Clear from proof of
previous proposition. $\blacksquare$

\begin{prop} Let $R$ be a finite rack and pick any $i$. Then, ${\cal O}\sb{i}$ is a subrack of $R$. \end{prop}Proof:
We will just note that ${\cal O}\sb{i}$ is closed under the $\ast$
and $\; \bar{\ast}\;$ operations. Consider the elements $(\cdots
((j\ast\sp{\epsilon\sb{1}}k\sb{1})\ast\sp{\epsilon\sb{2}}k\sb{2})\cdots
)\ast\sp{\epsilon\sb{m}}k\sb{m}$ and $(\cdots
((j'\ast\sp{\epsilon'\sb{1}}k'\sb{1})\ast\sp{\epsilon'\sb{2}}k'\sb{2})\cdots
)\ast\sp{\epsilon'\sb{m'}}k'\sb{m'}$ in ${\cal O}\sb{i}$. Then
\[
(\cdots
((j\ast\sp{\epsilon\sb{1}}k\sb{1})\ast\sp{\epsilon\sb{2}}k\sb{2})\cdots
)\ast\sp{\epsilon\sb{m}}k\sb{m}\ast\sp{\epsilon}(\cdots
((j'\ast\sp{\epsilon'\sb{1}}k'\sb{1})\ast\sp{\epsilon'\sb{2}}k'\sb{2})\cdots
)\ast\sp{\epsilon'\sb{m'}}k'\sb{m'}=(\cdots
((j\ast\sp{\epsilon''\sb{1}}k''\sb{1})\ast\sp{\epsilon''\sb{2}}k''\sb{2})\cdots
)\ast\sp{\epsilon''\sb{m''}}k''\sb{m''}
\]
since in the left-hand side of the equality, $j$ is the leftmost
element. Then ${\cal O}\sb{i}$ is closed under the $\ast$ and $\;
\bar{\ast}\;$ operations.
$\blacksquare$

\section{Indecomposability vs. Transitivity} \label{S:indecomptrans}

In this Section we briefly discuss indecomposability and
transitivity of racks and quandles.

\begin{def.} (\cite{AG}) A {\bf decomposition} of a rack $R=(X, \ast )$
is a disjoint union $X=X\sb 1 \sqcup X\sb 2$ such that $(X\sb 1,
\ast)$ and $(X\sb 2, \ast)$ are subracks of $R=(X, \ast )$ (in
particular, $X\sb 1$ and $X\sb 2$ are both non-empty). $R$ is {\bf
decomposable} if it admits a decomposition and {\bf
indecomposable} otherwise.
\end{def.}

\begin{def.} (\cite{Ryder}) A rack $R=(X, \ast )$ is said {\bf transitive}
if it is $c$-connected for some $c>0$.\end{def.}

\begin{prop} A finite rack is indecomposable if and only if it is transitive. \end{prop}Proof: See \cite{AG}. $\blacksquare$

Notice, moreover, that each of the two subracks in the definition
of decomposability above satisfy the following
\begin{def.}  A subrack $(X\sb o, \ast)$ of a rack $R=(X, \ast )$
is said {\bf stable} if, for any $i\sb o$ in $X\sb o$ and for any
$j$ in $X$, $i\sb o\ast j \in X\sb o$.\end{def.} (see \cite{AG}
for a proof). We remark, however, that a rack may be transitive
and still contain subracks. These subracks will not be stable, in
general. Consider, for instance, $R\sb 9$. It is a transitive rack
- see remarks after Proposition \ref{prop: 1-c}. On the other
hand, it contains three copies of $R\sb 3$ namely, $R\sp 1\sb 3 =
\{1, 4, 7 \}$, $R\sp 2\sb 3 = \{2, 5, 8\}$, $R\sp 3\sb 3 = \{3, 6,
9 \}$, see Table \ref{Ta:R9}.
\begin{table}[h!]
\begin{center}
    \begin{tabular}{| c || c | c | c | c | c | c | c | c | c |}\hline
              &   1     &      4  &     7   &    2   &   5   &   8   &    3   &   6   &   9   \\ \hline \hline
            1 &   {\bf 1}     &      {\bf 7}  &     {\bf 4}   &    3   &   9   &   6   &    5   &   2   &   8   \\ \hline
            4 &   {\bf 7}     &      {\bf 4}  &     {\bf 1}   &    9   &   6   &   3   &    2   &   8   &   5   \\ \hline
            7 &   {\bf 4}     &      {\bf 1}  &     {\bf 7}   &    6   &   3   &   9   &    8   &   5   &   2   \\ \hline
            2 &   9     &      6  &     3   &    {\bf 2}   &   {\bf 8}   &   {\bf 5}   &    4   &   1   &   7   \\ \hline
            5 &   6     &      3  &     9   &    {\bf 8}   &   {\bf 5}   &   {\bf 2}   &    1   &   7   &   4   \\ \hline
            8 &   3     &      9  &     6   &    {\bf 5}   &   {\bf 2}   &   {\bf 8}   &    7   &   4   &   1   \\ \hline
            3 &   8     &      5  &     2   &    1   &   7   &   4   &    {\bf 3}   &   {\bf 9}   &   {\bf 6}   \\ \hline
            6 &   5     &      2  &     8   &    7   &   4   &   1   &    {\bf 9}   &   {\bf 6}   &   {\bf 3}   \\ \hline
            9 &   2     &      8  &     5   &    4   &   1   &   7   &    {\bf 6}   &   {\bf 3}   &   {\bf 9}   \\ \hline
    \end{tabular}
\caption{$R\sb 9$ Multiplication Table}\label{Ta:R9}
\end{center}
\end{table}  Mind the unusual ordering of the underlying set and the boldface in order to bring out the three subracks. Moreover, these subracks are not stable: for instance, $1\ast 2= 3$ and $1 \in R\sp 1\sb 3$
whereas $3\notin R\sp 1\sb3$.

\section{Orbits and $\Delta$-Orbits} \label{S:orbs}

In this section we begin by describing an algorithm for obtaining
the orbits of a rack, given its sequence of permutations. This
algorithm is prompted by the following Proposition.

\begin{prop}\label{Prop:orbs-decomp} Let R be a finite rack of order $n$ with sequence of permutations $(\mu\sb 1,\dots ,\mu\sb n)$
 and assume the permutations are
already broken down into products of disjoint cycles. Let $\Sigma$
denote the set of all these disjoint cycles and suppose $\Sigma =
\Sigma\sb 1 \bigsqcup \Sigma\sb 2$ such that cycles from
$\Sigma\sb 1$ only move elements from $X\sb 1$ and cycles from
$\Sigma\sb 2$ only move elements from $X\sb 2$ where $X\sb 1$ and
$X\sb 2$ are disjoint subsets of $\{1, 2,\dots , n  \}$ the
underlying set of $R$. Then $R$ is decomposable.
\end{prop}Proof: Pick $x\sb i$ in $X\sb i$. Then, for any $j$ in $X$, $\mu\sb j(x\sb i)\in X\sb i$, for each
cycle in the permutations either moves an element of $X\sb 1$ or
(exclusively) an element of $X\sb 2$. $\blacksquare$

We now describe the algorithm for obtaining the orbits of a rack
given its sequence of permutations:

Fix a finite rack $R=(X, \ast)$ with sequence of permutations
$(\mu\sb 1,\dots ,\mu\sb n)$ and assume the permutations are
already broken down into products of disjoint cycles. Let $\Sigma
$ be the set of all cycles used in expressing the $\mu\sb i$'s as
products of disjoint cycles, except $1$-cycles. $\Sigma$ is a
finite set. If $\Sigma$ is the empty set, then each $\mu\sb{i}$ is
the product of $1$-cycles i.e.,
\[
\mu\sb{i}=(1)(2)(3)\cdots (n)=()
\]
and so $R$ is the trivial rack of order $n$, $T\sb{n}$. Otherwise,
assume $\Sigma = \{ \sigma\sb 1,\dots ,\sigma\sb m \}$ for some
positive integer $m$, where
\[
\sigma\sb k = (i\sb 1\sp k \cdots i\sb{m\sb k}\sp k) \qquad 1\leq k \leq m \qquad 1< m\sb k \leq n
\]
and $i\sb j\sp k \in \{ 1,\dots , n\}$. Set $s\sp k = \{ i\sb 1\sp
k, i\sb 2\sp k,\dots , i\sb{m\sb k}\sp k \}$, for each $k$. So,
the $s\sp k$ is the set of the elements which are moved by $\sigma
\sb k$. If $X \setminus \bigcup_{k=1}^{m}s\sp k$ is not empty, let
$\{ o\sb{1}, \dots , o\sb{l} \}=X \setminus \bigcup_{k=1}^{m}s\sp
k$. Then each $o\sb{i}$ is a (one element) orbit. Now for the
remaining orbits.

Let $E\sp 1\sb 1 = \{ i \: | \: s\sp 1 \cap s\sp i \neq \emptyset
\}$ and let $F\sp 1\sb 1 = \{ i \: | \: s\sp 1 \cap s\sp i =
\emptyset \}$.

Set
\[
\displaystyle {\cal O}\sb 1\sp 1 := \bigcup_{i\in E\sp 1\sb 1}s\sp
i
\]

If $ E\sp 1\sb 1 = \{ 1  \}$, then set ${\cal O}\sb 1={\cal O}\sb
1\sp 1$. Otherwise set $j\sp 1\sb 1 = \min \bigl( E\sp 1\sb 1
\setminus \{ 1 \}\bigr)$. Then let $E\sp 1\sb 2=\{ i\in F\sp 1\sb
1 \: | \: s\sp{j\sb 1}\cap s\sp i \neq \emptyset \}$, and $F\sp
1\sb 2=\{ i\in F\sp 1\sb 1 \: | \: s\sp{j\sp 1\sb 1}\cap s\sp i =
\emptyset \}$.

Set
\[
{\cal O}\sb 1\sp 2 := {\cal O}\sb 1\sp 1 \cup\bigcup_{i\in E\sp
1\sb 2}s\sp i
\]

If $ E\sp 1\sb 2 = \{ j\sp 1\sb 1  \}$, then set ${\cal O}\sb
1={\cal O}\sb 1\sp 2$. Otherwise set $j\sp 1\sb 2 = \min \bigl(
E\sp 1\sb 2 \setminus \{ j\sp 1\sb 1 \}\bigr)$. Then let $E\sp
1\sb 3=\{ i\in F\sp 1\sb 2 \: | \: s\sp{j\sp 1\sb 2}\cap s\sp i
\neq \emptyset \}$, and $F\sp 1\sb 3=\{ i\in F\sp 1\sb 2 \: | \:
s\sp{j\sp 1\sb 2}\cap s\sp i = \emptyset \}$.

Set
\[
{\cal O}\sb 1\sp 3 := {\cal O}\sb 1\sp 2 \cup\bigcup_{i\in E\sp
1\sb 3}s\sp i
\]

Since we are dealing with a finite rack, there has to be a
positive integer $k$ such that $ E\sp 1\sb k =\{ j\sb{k-1}\sp{1}
\}$. Then, set
\[
{\cal O}\sb 1= {\cal O}\sb 1\sp {k}
\]
and set $j\sp 2\sb 0 := \min F\sp 1\sb {k-1}$. This will be the
beginning of the construction of the second orbit.  Let $E\sp 2\sb
1 = \{ i \in F\sp 1\sb {k-1} \: | \: s\sp{j\sp 2\sb 0} \cap s\sp i
\neq \emptyset \}$ and let $F\sp 2\sb 1 = \{ i \in F\sp 1\sb {k-1}
\: | \: s\sp{j\sp 2\sb 0} \cap s\sp i = \emptyset \}$. Then set
\[
\displaystyle {\cal O}\sb 2\sp 1 := \bigcup_{i\in E\sp 2\sb 1}s\sp
i
\]

If $ E\sp 2\sb 1 = \{  j\sp 2\sb 0 \}$, then set ${\cal O}\sb
2={\cal O}\sb 2\sp 1$. Otherwise set $j\sp 2\sb 1 = \min \bigl(
E\sp 2\sb 1 \setminus \{ j\sp 2\sb 0 \}\bigr)$, and proceed as for
the previous orbit. Since we are dealing with a finite rack the
process has to stop at some step yielding, say, $N$ orbits, ${\cal
O}\sb 1,\dots , {\cal O}\sb N$ with more than one element.

The next results address the relations between profile and orbits
of racks.

\begin{prop} Let $R=(X, \ast )$ be a finite rack of order $n$. The
 permutations corresponding to the elements of a given orbit have
 all the same pattern.\end{prop}Proof: Let ${\cal O}\sb{i\sb 0}$ denote
the orbit of $i\sb 0$. Then for any $i \in {\cal O}\sb{i\sb 0}$,
we have $i =(\cdots ( (i\sb 0 \ast i\sb 1)\ast i\sb 2)\cdots )\ast
i\sb m$. Then
\begin{align}\notag
\mu\sb i &= \mu\sb{(\cdots ( (i\sb 0 \ast i\sb 1)\ast i\sb
2)\cdots )\ast i\sb m}=\mu\sb{\mu\sb{i\sb{m}}(\cdots ( (i\sb 0
\ast i\sb 1)\ast i\sb 2)\cdots )}=\mu\sb{i\sb{m}}\mu\sb{((\cdots (
(i\sb 0 \ast i\sb 1)\ast i\sb 2)\cdots
))}\mu\sb{i\sb{m}}\sp{-1}=\dots =\\ \notag &=\mu\sb{i\sb m}\cdots
\mu\sb{i\sb 1}\mu\sb{i\sb 0}\mu\sb{i\sb 1}\sp{-1}\cdots
\mu\sb{i\sb m}\sp{-1}
\end{align}
Since the permutation corresponding to any $i$ in ${\cal
O}\sb{i\sb 0}$ is conjugate to the permutation corresponding to
$i\sb{0}$ then they have the same pattern. The result follows.
$\blacksquare$

\begin{cor} Let $R=(X, \ast )$ be an indecomposable finite rack of order $n$. Then its permutations have all
 the same pattern. In particular, the profile is a constant sequence and the detail is a singleton.
 \end{cor}Proof: Since an indecomposable rack is transitive, the result follows from the proof of the previous
 proposition. $\blacksquare$

\begin{def.}\label{def:deltaorbit} Let  $R=(X, \ast )$ be a finite rack of order $n$. For each $i$ in $\{1, 2,\dots
,n\}$ consider the sequence $(i\sb k)$ defined in the following
way. $i\sb 1= i$, and, for any $k\geq 1$, $i\sb{k+1}=i\sb k\ast
i\sb k$. We call the set of its terms the {\bf $\Delta$-orbit of
$i$}, or simply the {\bf $\Delta\sb i$ orbit}. The length of the
$\Delta$-orbit is the number of elements in it.
\end{def.}

We remark that the $\Delta$-orbits can be defined for infinite
racks, but, if we are dealing with finite racks, then the
$\Delta\sb i$ orbits are finite sets. Furthermore, for quandles,
$\Delta$-orbits are singletons.

\begin{prop}\label{prop:periodic-deltaorbit} Let $R=(X, \ast )$ be a finite rack of order $n$. For an arbitrary $i$
in $\{1,2, \dots , n \}$, consider the sequence introduced in
Definition \ref{def:deltaorbit}. This sequence is periodic and the
length of $\Delta\sb{i}$ is its period. Moreover, for any $k\geq
1$, $i\sb{k+1}=\mu\sb{i}\sp{k}(i)$.\end{prop}Proof: We start out
by proving the last claim, for any $k\geq 1$,
$i\sb{k+1}=\mu\sb{i}\sp{k}(i)$. We have $i\sb{2}=i\ast i =
\mu\sb{i}(i)=\mu\sb{i}\sp{1}(i)$. Suppose there is a positive
integer $r$ such that $i\sb{r+1}=\mu\sb{i}\sp{r}(i)$. Then,
\begin{align}\notag
i\sb{r+2}&=i\sb{r+1}\ast
i\sb{r+1}=\mu\sb{i\sb{r+1}}(i\sb{r+1})=\mu\sb{\mu\sb{i}\sp{r}(i)}(\mu\sb{i}\sp{r}(i))=
\mu\sb{\mu\sb{i}(\mu\sb{i}\sp{r-1}(i))}(\mu\sb{i}\sp{r}(i))=\\
\notag
&=\mu\sb{i}\mu\sb{\mu\sb{i}\sp{r-1}(i)}\mu\sb{i}\sp{-1}(\mu\sb{i}\sp{r}(i))=\dots
=\mu\sb{i}\sp{r}\mu\sb{i}\mu\sb{i}\sp{-r}(\mu\sb{i}\sp{r}(i))=\mu\sb{i}\sp{r+1}(i)
\end{align}
thus proving the claim.

Since $R$ is finite, there is a positive integer $m$ such that
$i\sb{m}\ast i\sb{m}$ is in $\{ i\sb{k}\; |\; 1\leq k \leq m \}$.
Suppose $i\sb{m}\ast i\sb{m}=i\sb{l}$ with $l$ not equal to $1$.
Then $i\sb{l+1}=i\sb{l}\ast i\sb{l},\;  \dots \; ,
i\sb{m}=i\sb{m-1}\ast i\sb{m-1}, \; i\sb{l}=i\sb{m}\ast i\sb{m}\;
\dots \; $, and so on, keeping on repeating the subsequence $(
i\sb{l},\dots , i\sb{m})$. So for any $k\geq m-1$,
\[
i=i\sb{1}\neq i\sb{k+1}\ast
i\sb{k+1}=\mu\sb{i\sb{k+1}}(i\sb{k+1})=\mu\sb{\mu\sb{i}\sp{k}(i)}(\mu\sb{i}\sp{k}(i))
=\mu\sb{i}\sp{k}\mu\sb{i}\mu\sb{i}\sp{-k}(\mu\sb{i}\sp{k}(i))=\mu\sb{i}\sp{k+1}(i)
\]
Via Proposition \ref{prop:musp{n}1} we know that there are
infinitely many positive integers $r$ such that $\mu\sb{i}\sp{r}$
is the identity permutation which clearly contradicts the
preceding calculations. Hence $l=1$ i.e., the sequence $i=i\sb{1}$
and $i\sb{k+1}=i\sb{k}\ast i\sb{k}$ for any $k\geq 1$, is periodic
and the length of the $\Delta\sb{i}$ orbit is its period.
$\blacksquare$

\begin{prop}\label{prop:perm-deltaorbit} Let $R=(X, \ast )$ be a finite rack of order $n$.
The permutations corresponding to the elements of a given
$\Delta$-orbit are all equal. Furthermore, for any $i$ in $\{1,
2,\dots , n \}$, $\mu\sb{i}=(i\sb{1} \; i\sb{2}\dots i\sb{m})\pi$,
where $m$ is the length of the $\Delta$-orbit, $(i\sb{1}, i\sb{2},
\dots , i\sb{m})$ is the sequence introduced in Definition
\ref{def:deltaorbit}, and $\pi$ is a permutation of the remaining
elements i.e., $\{1, 2,\dots , n \}\setminus \{i\sb{1},
i\sb{2},\dots , i\sb{m} \}$
\end{prop}Proof: Fix an arbitrary $i$ in $\{1, 2,\dots , n \}$
 and consider the sequence $(i\sb k)$ introduced in Definition \ref{def:deltaorbit}.
 We have $\mu\sb{i\sb 2}=\mu\sb{i\ast i}=\mu\sb{\mu\sb i(i)}=\mu\sb i\mu\sb i\mu\sp{-1}\sb i = \mu\sb i$.
  Suppose $\mu\sb{i\sb k}=\mu\sb i$. Then
  $\mu\sb{i\sb{k+1}}=\mu\sb{i\sb k\ast i\sb k}=\mu\sb{\mu\sb{i\sb k}(i\sb k)}=\mu\sb{i\sb k}\mu\sb{i\sb k}\mu\sp{-1}\sb{i\sb k}=\mu\sb{i\sb k}=\mu\sb
  i$ which establishes the first claim in the statement. For the second claim, we know from the previous Proposition
  that $i\sb{k+1}=\mu\sb{i}\sp{k}(i)$, for all $k \geq 1$. Then,
  either $m=1$ and $\mu\sb{i}(i)=i$ or $m>1$ and $i$ is moved to
  $\mu\sb{i}(i)$ which is moved to $\mu\sb{i}\sp{2}(i)$ and so on
  and so forth until $i$ is reached again. The proof is complete.
   $\blacksquare$

\begin{prop}\label{prop:deltaorbits:disjoint} Let $R=(X, \ast )$ be a finite rack of order $n$.
 Suppose there are elements $i$ and $j$ in $X$ such that
 $\Delta\sb i \cap \Delta\sb j$ is non-empty. Then $\Delta\sb i = \Delta\sb j$.
 \end{prop}Proof: Let $\Delta\sb i
= \{ i\sb m | m\in \mathbb{Z}\sp{+} \}$ and $\Delta\sb j
 = \{ j\sb m | m\in \mathbb{Z}\sp{+} \}$ where the $i\sb k$ and $j\sb k$
 have the same meaning as before. Since $R$ is a finite rack,
 these sequences are periodic. Since their intersection is not
 empty then there exist $k$ and $k'$ in $\mathbb{Z}\sp{+}$ such that $i\sb k = j\sb
 {k'}$. Then, by induction, if there exists $p$ such that
 $i\sb{k+p}=j\sb{k'+p}$, then $i\sb{k+p}\ast i\sb{k+p}=j\sb{k'+p}\ast
 j\sb{k'+p}$ i.e., $i\sb{k+p+1}=j\sb{k'+p+1}$ and since the $(i\sb
 m)$ and $(j\sb m)$ sequences are periodic, the result follows.   $\blacksquare$

\begin{prop} \label{prop:deltaorbit:iso}Let $\alpha$ be a rack isomorphism between
racks $R$ and $R'$. For any $i$, let $\Delta\sb{i}$ and
$\Delta'\sb{i}$ denote the $\Delta$-orbit of $i$ in $R$ and $R'$,
respectively. Then $\alpha$ restricted to $\Delta\sb{i}$ is a
bijection of $\Delta\sb{i}$ with $\Delta'\sb{\alpha (i)}$.
Furthermore, $\alpha$ induces a bijection from the set of
$\Delta$-orbits of $R$ to the set of $\Delta$-orbits of $R'$.
\end{prop}Proof: Let $\Delta\sb{i}=\{ i\sb{k} \; | \; k\in \mathbb{Z}\sp{+} \}$,
where the sequence $(i\sb{k})$ has the same meaning as in Definition
\ref{def:deltaorbit}. Then
\[
 \alpha \bigl( \Delta\sb{i}\bigr)  =\alpha \Big( \{ i\sb{k} \; |
\; k\in \mathbb{Z}\sp{+} \}\Bigr)=
 \{ \alpha (i\sb{k}) \; | \; k\in \mathbb{Z}\sp{+}
 \}=\Delta\sb{\alpha (i)}
\]
noting that $\alpha (i\sb{1})=\alpha (i)$ and, for any $k$,
 $\alpha (i\sb{k+1})=\alpha (i\sb{k}\ast i\sb{k})=\alpha (i\sb{k})\ast' \alpha (i\sb{k})$. $\blacksquare$

\begin{prop}\label{prop:deltaorbit:subrack} Let $R$ be a finite rack of order $n$ and pick
$i$ in $\{ 1, 2,\dots , n \}$. Then $\Delta\sb{i}$, the $\Delta$-orbit of $i$,
is a subrack of $R$.\end{prop}Proof: We will just show that
$\Delta\sb{i}$ is closed under the $\ast$ operation. Pick any two
elements in $\Delta\sb{i}$, $i\sb{k'}$ and $i\sb{k}$. Due to the
periodicity of the sequence $(i\sb{n})$ then $i\sb{k'}=i\sb{k+p}$,
for a given $p$. We will prove by induction on $p$ that
$i\sb{k'}\ast i\sb{k}(=i\sb{k+p}\ast i\sb{k})=i\sb{k+p}\ast
i\sb{k+p}$, thereby proving it is an element of $\Delta\sb{i}$.
For $p=1$ we have, $i\sb{k+1}\ast i\sb{k}=\big( i\sb{k}\ast
i\sb{k} \big)\ast i\sb{k}=\big( i\sb{k}\ast i\sb{k} \big)\ast
\big( i\sb{k}\ast i\sb{k} \big)=i\sb{k+1}\ast i\sb{k+1}$. Now
assume that for some $p$ it is true that $i\sb{k+p}\ast
i\sb{k}=i\sb{k+p}\ast i\sb{k+p}$. Then, $i\sb{k+p+1}\ast
i\sb{k}=\big( i\sb{k+p}\ast i\sb {k+p}\big)\ast i\sb{k}=\big(
i\sb{k+p}\ast i\sb {k}\big)\ast\big( i\sb{k+p}\ast i\sb
{k}\big)=\big( i\sb{k+p}\ast i\sb {k+p}\big)\ast\big(
i\sb{k+p}\ast i\sb {k+p}\big)=i\sb{k+p+1}\ast i\sb {k+p+1}$, an
element of $\Delta\sb{i}$.
 $\blacksquare$

\section{Racks Given Their Profiles} \label{S:racksfromprof}

In this section we express the possible representatives of the
equivalence classes of racks given specific profiles. As usual,
let $R=(X, \ast )$ be a rack of order $n$ and $(\mu\sb 1,\dots
,\mu\sb n)$ be its sequence of permutations.

\begin{thm} Suppose one of the permutations of a finite rack of order $n$ is a cycle of length $n$.
Then, the sequence of permutations of this rack is constant.
Moreover there is only one such rack of this order (up to
isomorphism). This rack can be taken to be the $C\sb n$ rack i.e.,
the cyclic rack of order $n$.\end{thm}Proof: Let $\mu\sb i=(i\sb
1\;  i\sb 2 \; \cdots \; i\sb n)$ where, without loss of
 generality, we take $i\sb 1 = i$. Then, by Proposition \ref{prop:perm-deltaorbit}, $\{ 1, 2,\dots , n
 \}$ is the $\Delta\sb{i}$-orbit which contains all elements of
 the rack. Again by Proposition \ref{prop:perm-deltaorbit}, $\mu\sb j=(i\sb
 1\; i\sb 2 \; \cdots \; i\sb n)$, for all $j$ in $\{ 1, 2,\dots , n \}$.

Now consider the assignment $\alpha(i\sb{k})=k$, for all $k$ in
$\{ 1, 2,\dots , n \}$. This is a rack isomorphism from the rack
under study to the cyclic rack of order $n$, $C\sb n$, which
finishes the proof. $\blacksquare$

We next address the case of a prack of order $n$ with constant
profile $(\{m, n-m\},\dots , \{m, n-m\})$. For clarity, we split
this study into Theorem \ref{th:prack1,n-1}, concerning $m=1$,
Corollary \ref{cor:prackm,n-m}, concerning $1<m<n-m$, and
Corollary \ref{cor:prackm,m}, concerning $m=n-m$. We recall that
the Euler $\varphi$-function is, for each positive integer $n$,
the number of positive integers relatively prime to, and less than
$n$.

\begin{thm}\label{th:prack1,n-1} Any prack of order $n>2$ with constant profile
\[
(\{1, n-1\},\dots , \{1, n-1\})
\]
is decomposable and is isomorphic to a prack given by the
following sequence of permutations:
\begin{align}\notag
\mu\sb{1}&=\mu\sb{2}=\quad \cdots \quad =\mu\sb{n-1}=( 1\; 2\;
\cdots \; n-1 )(n),\\ \notag \mu\sb{n}&=( 1\; 2\; \cdots \; n-1
)\sp{k}(n)
\end{align}
where $k$ is relatively prime to $n-1$.

Moreover, there are exactly $\varphi (n-1)$ such pracks up to
isomorphism, one per each $k$ relatively prime to $n-1$ ($1\leq k
< n-1$).\end{thm}Proof: Since each permutation fixes one element,
then there has to be one which does not fix its own index, for
otherwise the rack would be a quandle. So suppose
$\mu\sb{i}(j)=j$. If necessary using the rack isomorphism $\alpha
= (i \; 1)(j\; n)$ we may assume $\mu\sb 1 (n)=n$. Let $i\sb{1}=1$
and let the sequence $(i\sb{k})$ be as in Definition
\ref{def:deltaorbit}. Then by Proposition
\ref{prop:perm-deltaorbit}, $\mu\sb{1}=\mu\sb{2}=\cdots
=\mu\sb{n-1}=(i\sb{1} \; i\sb{2}\; \cdots \; i\sb{n-1})(n)$.
Moreover, $\mu\sb n(n) = n$ for otherwise suppose $\mu\sb n(n) =
i\sb n \neq n$. Then $\mu\sb 1 = \mu\sb{i\sb n}=\mu\sb{\mu\sb
n(n)} = \mu\sb n\mu\sb n\mu\sb n\sp{-1} = \mu\sb n$ and so $\mu\sb
n (n)=\mu\sb 1(n)= n$. In particular, this prack is decomposable
by Proposition \ref{Prop:orbs-decomp}. Using, if necessary, the
isomorphism $\alpha (i\sb{k})=k$, for $1\leq k \leq n-1$, the
first $n-1$ permutations of the sequence of permutations become
all equal to $\mu\sb{k}= ( 1\;  2 \; \cdots \; n-1 )( n )$, for
$1\leq k \leq n-1$. For this prack the conditions
$\mu\sb{\mu\sb{j}(i)}=\mu\sb j\mu\sb i\mu\sp{-1}\sb j$ are
equivalent to saying that $\mu\sb{n}$ commutes with $( 1\; 2 \;
\cdots \; n-1 )( n )$ which is equivalent to saying that the
$(n-1)$-cycle in $\mu\sb{n}$ is in the centralizer of $( 1\;  2\;
\cdots \; n-1 )$ in $S\sb{n-1}$. The centralizer of an $r$-cycle
in $S\sb{r}$ is the cyclic group generated by that cycle. Thus, in
the current situation we want to pick a power $( 1\;  2 \dots n-1
)\sp{k}$ which is an $(n-1)$-cycle which amounts to picking the
$k$ relatively prime to $n-1$. Hence the prack under study can
only be isomorphic to a prack given by a sequence of permutations
as in the statement and there will be at most $\varphi (n-1)$
pracks up to isomorphism with the indicated profile. We will next
show that there are exactly $\varphi (n-1)$ such pracks by showing
that any two such pracks with different $k$'s are not isomorphic.

Suppose two such pracks are isomorphic. Specifically, assume that
rack $R$ with permutations $\mu\sb{1}=\dots =\mu\sb{n-1}=(1\;
 2 \; \cdots \; n-1 )(n)$, $\mu\sb{n}=(1\;  2 \; \cdots \; n-1 )\sp{k}(n)$ is
isomorphic with rack $R'$ with permutations $\mu'\sb{1}=\dots
=\mu'\sb{n-1}=(1\;  2 \; \cdots \; n-1 )(n)$, $\mu'\sb{n}=(1\;  2
\; \cdots \; n-1 )\sp{k'}(n)$, via an isomorphism $\alpha $ (with
distinct integers $1\leq k, k'< n-1$ both relatively prime to
$n-1$). Since, by Proposition \ref{prop:deltaorbit:iso},
isomorphisms map $\Delta$-orbits to $\Delta$-orbits of the same
length, then $\alpha (n)=n$, and $\alpha (\{1, 2, \dots ,
n-1\})=\{1, 2, \dots , n-1\}$. Hence, for any $i$ in $\{1, 2,
\dots , n-1\}$,
\[
(1\;  2 \; \cdots \; n-1 )(n)=\mu'\sb{\alpha(i)}=\alpha
\mu\sb{i}\alpha\sp{-1}=\alpha (1\;  2 \; \cdots \; n-1
)(n)\alpha\sp{-1}
\]
i.e., $\alpha $ is in the centralizer of $(1\;  2 \; \cdots \; n-1
)(n)$ in $S\sb{n}$. Since $\alpha = \sigma \cdot (n)$ for some
$\sigma$ in $S\sb{n-1}$ then $\sigma$ is in the centralizer of
$(1\;  2 \; \cdots \; n-1 )$ in $S\sb{n-1}$. Thus, $\sigma = ( 1\;
2 \; \cdots \; n-1 )\sp{l}$ for some $l$ in $\{1, 2, \dots ,
n-1\}$. Finally,
\begin{align}\notag
(1\;  2 \; \cdots \; n-1 )\sp{k'}(n)&=\mu'\sb{n}=\mu'\sb{\alpha
(n)}=\alpha \mu\sb{n}\alpha\sp{-1}=(1\;  2 \; \cdots \; n-1
)\sp{l}(n)(1 \; 2 \; \cdots \; n-1 )\sp{k}(n)(1\;  2 \; \cdots \;
n-1 )\sp{-l}(n)=\\ \notag
 &=(1\;  2 \; \cdots \; n-1)\sp{k}(n)
\end{align}
Hence $k=k'$ and the result follows.
$\blacksquare$

We remark there are indecomposable quandles of order $n$ of
constant profile $\{ (1, n-1),\dots , (1, n-1) \}$ - see, for
instance, the quandles $\mathbb{Z}\sb{5}[T, T\sp{-1}]/(T-2)$ and
$\mathbb{Z}\sb{5}[T, T\sp{-1}]/(T-3)$ in Section \ref{S:inv2}.

\begin{cor}\label{cor:prackm,n-m} Any prack of order $n>3$ with constant profile
\[
(\{m, n-m\},\dots , \{m, n-m\})
\]
(where $1<m< n-m$) is decomposable and is isomorphic to a prack
given by the following sequence of permutations:
\begin{align}\notag
\mu\sb{1}&=\mu\sb{2}=\quad  \dots \quad =\mu\sb{m}=( 1\; 2\;
\cdots \; m )( m+1\; \cdots \; n )\sp{l}, \\ \notag
\mu\sb{m+1}&=\mu\sb{m+2}=\quad \dots \quad =\mu\sb{n}=( 1\; 2\;
\cdots \; m )\sp{k}( m+1\; \cdots \; n )
\end{align}
where $k$ is relatively prime to $m$ and $l$ is relatively prime
to $n-m$.

There are exactly $\varphi (m)\cdot \varphi (n-m)$ such pracks up
to isomorphism, one prack per each ordered pair $(k, l)$ where $k$
is relatively prime to $m$ and $l$ is relatively prime to $n-m$
($1\leq k < m$, $1\leq l < n-m$).\end{cor}Proof: $\mu\sb{1}(1)\neq
1$ for there is no cycle of length one. Let $i\sb{1}=1$ and let
the sequence $(i\sb{k})$ be as in Definition \ref{def:deltaorbit}.
Then by Proposition \ref{prop:perm-deltaorbit},
$\mu\sb{i\sb{1}}=\mu\sb{i\sb{2}}=\cdots =\mu\sb{i\sb{m}}=(i\sb{1}
\; i\sb{2}\; \cdots \; i\sb{m})\pi$, where $\pi$ is a permutation
of $\{1, 2,\dots , n \}\setminus \{i\sb{1}, i\sb{2},\dots ,
i\sb{m} \}$. Pick any $j$ in $\{1, 2,\dots , n \}\setminus
\{i\sb{1}, i\sb{2},\dots , i\sb{m} \}$. Again, $\mu\sb{j}(j)\neq
j$ for there is no cycle of length one. Let $j\sb{1}=j$ and
$(j\sb{k})$ be as in Definition \ref{def:deltaorbit}. Arguing as
above, $\mu\sb{j\sb{1}}=\mu\sb{j\sb{2}}=\cdots
=\mu\sb{j\sb{m}}=(j\sb{1} \; j\sb{2}\; \cdots \; j\sb{n-m})\rho$,
where $\rho$ is a permutation of $\{i\sb{1}, i\sb{2},\dots ,
i\sb{m} \}$. Let $\alpha$ be the isomorphism such that
$\alpha(i\sb{k})=k$ and $\alpha(j\sb{k})=m+k$. In this way the
permutations are now (keeping the notation $\mu$ for the
permutations in the image of the isomorphism)
$\mu\sb{1}=\mu\sb{2}=\cdots =\mu\sb{m}=(1 \; 2\; \cdots \;
m)\delta$ and $\mu\sb{m+1}=\mu\sb{m+2}=\cdots =\mu\sb{n}=\gamma
(m+1 \; m+2\; \cdots \; n)$, where $\delta$ is a permutation of
$\{m+1, m+2,\dots , n \}$ and $\gamma$ is a permutation of $\{1,
2,\dots , m\}$. The self-distributivity axiom amounts to saying
that $\mu\sb{\mu\sb{i}(j)}=\mu\sb{i}\mu\sb{j}\mu\sb{i}\sp{-1}$. In
this way, for $i$, $j$ in $\{1, 2,\dots , m\}$, or $i$, $j$ in
$\{m+1, m+2,\dots , n \}$ they yield trivial expressions. If $i$
is from $\{1, 2,\dots , m\}$ and $j$ is from $\{m+1, m+2,\dots , n
\}$ we obtain, since $\delta (j)$ is in $\{m+1, m+2,\dots , n \}$,
\begin{align}\notag
\gamma (m+1 \; m+2\; \dots \; n)=\mu\sb{\delta(j)}&=
\mu\sb{\mu(i)}=(1 \; 2\; \dots \; m)\delta \cdot \gamma (m+1 \;
m+2\; \cdots \; n)\cdot \delta\sp{-1}(m \; \cdots \; 2 \; 1)=
\\ \notag &=(1 \; 2\; \cdots \; m)\gamma (m \; \cdots \; 2 \;
1)\cdot \delta (m+1 \; m+2\; \cdots \; n)\cdot \delta\sp{-1}\notag
\end{align}
so $\gamma =(1 \; 2 \; \dots \; m)\gamma (m \; \cdots \; 2 \; 1)$
and $(m+1 \; m+2 \; \dots \; n)=\delta (m+1 \; m+2 \; \cdots \;
n)\cdot \delta\sp{-1}$ i.e., $\gamma $ is in the centralizer of
$(1 \; 2 \; \dots \; m)$ in the symmetry group of $\{1, 2,\dots ,
m\}$, and $\delta $ is in the centralizer of $(m+1 \; m+2 \;
\cdots \; n)$ in the symmetry group of $\{m+1, m+2,\dots , n \}$.
Thus, $\gamma =(1 \; 2 \; \dots \; m)\sp{k}$, for some $k$
relatively prime to $m$, and $\delta =(m+1 \; m+2 \; \cdots \;
n)\sp{l}$, for some $l$ relatively prime to $n-m$, since $\gamma$
(resp., $\delta$) has to be an $m$-cycle (resp., $(n-m)$-cycle).
We arrive at the same conclusion if $i$ is from $\{m+1, m+2,\dots
, n \}$ and $j$ is from $\{1, 2,\dots , m\}$. In particular, this
rack is indecomposable by Proposition \ref{Prop:orbs-decomp}.

Now suppose two such racks are isomorphic. Specifically, assume
that rack $R$ with permutations $\mu\sb{1}=\dots =\mu\sb{m}=(1\;
 2 \; \cdots \; m )(m+1\;
 2 \; \cdots \; n)\sp{l}$, $\mu\sb{m+1}=\dots =\mu\sb{n}=(1\;
 2 \; \cdots \; m )\sp{k}(m+1\;
 2 \; \cdots \; n)$ is
isomorphic with rack $R'$ with permutations $\mu'\sb{1}=\dots
=\mu'\sb{m}=(1\; 2 \; \cdots \; m )(m+1\;
 2 \; \cdots \; n)\sp{l'}$, $\mu'\sb{m+1}=\dots =\mu'\sb{n}=(1\;
 2 \; \cdots \; m )\sp{k'}(m+1\;
 2 \; \cdots \; n)$, via an isomorphism $\alpha $ (with  distinct pairs
 $( k, l)$, and $( k', l')$ further satisfying the conditions in
 the statement).
Since, by Proposition \ref{prop:deltaorbit:iso}, isomorphisms map
$\Delta$-orbits to $\Delta$-orbits of the same length, then
$\alpha (\{1, 2, \dots , m \})=\{1, 2, \dots , m \}$, and $\alpha
(\{m+1, \dots , n \})=\{m+1, \dots , n \}$. Hence, $\alpha =
\delta \gamma$, where $\delta$ is in $S\sb{\{1, 2, \dots , m \}}$
and $\gamma$ is in $S\sb{\{m+1, \dots , n \}}$. Then, for each $i$
in $\{1, 2, \dots , m \}$, there is an $i'$ in $\{1, 2, \dots , m
\}$ such that $\alpha (i)=i'$; and, for each $j$ in $\{m+1, \dots
, n \}$, there is a $j'$ in $\{m+1, \dots , n \}$ such that
$\alpha (j)=j'$. So
\[
(1\; 2 \; \cdots \; m )(m+1\; 2 \; \cdots \;
n)\sp{l'}=\mu'\sb{\alpha(i)}=\alpha \mu\sb{i}\alpha\sp{-1}=\delta
\gamma (1\; 2 \; \cdots \; m )(m+1\;  \cdots \; n)\sp{l}
\delta\sp{-1} \gamma\sp{-1}
\]
and
\[
(1\; 2 \; \cdots \; m )\sp{k'}(m+1\; 2 \; \cdots \;
n)=\mu'\sb{\alpha(j)}=\alpha \mu\sb{j}\alpha\sp{-1}=\delta \gamma
(1\; 2 \; \cdots \; m )\sp{k}(m+1\;  \cdots \; n) \delta\sp{-1}
\gamma\sp{-1}
\]
whence $\delta $ is in the centralizer of $(1\;  2 \; \cdots \; m
)$ in $S\sb{\{1, 2, \dots , m \}}$, $\gamma $ is in the
centralizer of $(m+1\; \cdots \; n )$ in $S\sb{\{m+1, \dots , n
\}}$; and $l=l'$, and $k=k'$ i.e, the racks are isomorphic. In
particular, there are $\varphi (m)\cdot \varphi (n-m)$ pracks up
to isomorphism with the indicated profile. This concludes the
proof. $\blacksquare$

\begin{cor}\label{cor:prackm,m} Any prack of order $n>3$ with constant profile
\[
(\{m, m\},\dots , \{m, m\})
\]
(where $m>1$ and $n=m+m$) is decomposable and is isomorphic to a
prack given by the following sequence of permutations:
\begin{align}\notag
\mu\sb{1}&=\mu\sb{2}=\quad  \dots \quad =\mu\sb{m}=( 1\; 2\;
\cdots \; m )( m+1\; \cdots \; n )\sp{l}, \\ \notag
\mu\sb{m+1}&=\mu\sb{m+2}=\quad \dots \quad =\mu\sb{n}=( 1\; 2\;
\cdots \; m )\sp{k}( m+1\; \cdots \; n )
\end{align}
where $k$ and $l$ are both relatively prime to $m$.

There are exactly $\frac{1}{2}(1+\varphi (m))\cdot \varphi (m)$
such pracks up to isomorphism, one per each ordered pair $(k, l)$
where $k$ and $l$ are both relatively prime to $m$, and $1\leq l
\leq k<m$.\end{cor}Proof: We omit the proof here, since it is a
slight variation to the proof to the previous Corollary. We just
note that, since the $\Delta$-orbits have the same length, then
different choices of $k$ and $l$ may give rise to isomorphic
pracks. This is taken care of, in the statement, by requiring
$k\geq l$. $\blacksquare$

We now address the case of pracks of constant profile such that
the corresponding permutations have more than two cycles. Theorem
\ref{th:prack1,1,...,n-f} and Corollary \ref{cor:prackf,n} concern
pracks of constant profile with $r$ cycles of length one and one
cycle of length greater than one. Theorem \ref{th:prackm1<m2<m3}
concerns a prack of constant profile with three cycles of distinct
lengths, all greater than one.

\begin{thm}\label{th:prack1,1,...,n-f}Given integers $n>3$
and $f>2$ such that $f<n$, the sequence of permutations,
\begin{align}\notag
\mu\sb{1}&=\mu\sb{2}=\dots =\mu\sb{n-f}=(1 \; 2 \; \cdots \;
n-f)(n-f+1)(n-f+2)\cdots (n-1)(n)\\ \notag \mu\sb{j}&=(1 \; 2 \;
\dots \; n-f)\sp{k\sb{j}}(n-f+1)(n-f+2)\cdots (n-1)(n), \quad j\in
\{ n-f+1,\dots , n \}
\end{align}
where the $k\sb{j}$'s are relatively prime to $n-f$, stands for a
decomposable prack of order $n$ with profile:
\[
(\{1, 1, \dots , 1, n-f \},\dots , \{1, 1, \dots , 1, n-f \})
\]

Let $P\sb{i}(f)$ be the number of partitions of $f$ into $i$
summands, $C\sb{i}\sp{\varphi (n-f)}$ be the number of subsets of
$i$ elements from a set of $\varphi (n-f)$ elements, and $m:=\min
\{ f, \varphi (n-f) \}$. There are exactly
\[
\sum\sb{i=1}\sp{m}C\sb{i}\sp{\varphi (n-f)}P\sb{i}(f)
\]
pracks with such sequence of permutations up to isomorphism, one
per finite sequence $\{ k\sb{j} \}\sb{j=n-f+1}\sp{n}$, modulo
permutation of the $k\sb{j}$'s ($k\sb{j}$'s relatively prime to,
and less than, $n-f$).\end{thm}

Before proving the Theorem we state:

\begin{cor}\label{cor:prackf,n} Consider the profile in the statement of Theorem
\ref{th:prack1,1,...,n-f}:
\[
(\{1, 1, \dots , 1, n-f \},\dots , \{1, 1, \dots , 1, n-f \})
\]

If $2<f<n-f$ there are exactly
$\sum\sb{i=1}\sp{m}C\sb{i}\sp{\varphi (n-f)}P\sb{i}(f)$ pracks of
order $n$ with such a profile, up to isomorphism. These pracks are
the ones indicated in the Theorem.

If $f=n-f$, there are exactly
$1+\sum\sb{i=1}\sp{m}C\sb{i}\sp{\varphi (n-f)}P\sb{i}(f)$ pracks
of order $n$ with such a profile, up to isomorphism. These pracks
are the ones indicated in the Theorem plus the decomposable prack
given by the following sequence of permutations:
\begin{align}\notag
\mu\sb{1}&=\mu\sb{2}=\cdots =\mu\sb{n-f}=(1 \; 2 \; \cdots \;
n-f)(n-f+1)(n-f+2)\cdots (n-1)(n)\\ \notag
\mu\sb{n-f+1}&=\mu\sb{n-f+1}=\cdots =\mu\sb{n}=(1)(2)\cdots
(n-f)(n-f+1 \; n-f+2 \; \cdots \; n-1 \; n)
\end{align}\end{cor}

Proof (Theorem \ref{th:prack1,1,...,n-f}): Since we are dealing
with a prack, there has to be a permutation which moves its own
index. It follows that the permutations corresponding to the
$\Delta$-orbit of this element are all equal. Thus without loss of
generality we may assume, arguing as in similar cases above,
\[
\mu\sb{1}=\mu\sb{2}=\cdots =\mu\sb{n-f}=(1 \; 2 \; \cdots \;
n-f)(n-f+1)(n-f+2)\cdots (n-1)(n)
\]
We now assume that the remaining permutations fix the same
elements as the former ones do. In this way,
\[
\mu\sb{j}=\sigma\sb{j}(n-f+1)(n-f+2)\cdots (n-1)(n)
\]
for $j$ in $\{ n-f+1,\dots , n \}$, where $\sigma\sb{j}$ is an
$(n-f)$-cycle in $S\sb{\{ 1,\dots , n-f \}}$. Using the
distributivity axiom, $\sigma\sb{j}$ has to be in the centralizer
of $(1\; 2\; \cdots \; n-f)$ in $S\sb{ \{ 1,\dots , n-f \} }$.
Hence,
\[
\mu\sb{j}=(1 \; 2 \; \cdots \;
n-f)\sp{k\sb{j}}(n-f+1)(n-f+2)\cdots (n-1)(n)
\]
with $k\sb{j}$ relatively prime to $n-f$, for each $j\in \{
n-f+1,\dots n \}$. Clearly, this prack is decomposable by
Proposition \ref{Prop:orbs-decomp}.

Suppose two such racks are isomorphic. Specifically, assume that
\begin{align}\notag
\mu\sb{1}&=\mu\sb{2}=\dots =\mu\sb{n-f}=(1 \; 2 \; \cdots \;
n-f)(n-f+1)(n-f+2)\cdots (n-1)(n)\\ \notag \mu\sb{j}&=(1 \; 2 \;
\dots \; n-f)\sp{k\sb{j}}(n-f+1)(n-f+2)\cdots (n-1)(n), \quad j\in
\{ n-f+1,\dots , n \}
\end{align}
and
\begin{align}\notag
\mu'\sb{1}&=\mu'\sb{2}=\dots =\mu'\sb{n-f}=(1 \; 2 \; \cdots \;
n-f)(n-f+1)(n-f+2)\cdots (n-1)(n)\\ \notag \mu'\sb{j}&=(1 \; 2 \;
\dots \; n-f)\sp{k'\sb{j}}(n-f+1)(n-f+2)\cdots (n-1)(n), \quad
j\in \{ n-f+1,\dots , n \}
\end{align}
(where $\{ k\sb{j} \}\sb{j=n-f+1}\sp{n}$ and $\{ k'\sb{j}
\}\sb{j=n-f+1}\sp{n}$ are sequences of numbers relatively prime
to, and less than $n-f$) are isomorphic via isomorphism $\alpha$.
Since rack isomorphisms map $\Delta$-orbits into $\Delta$-orbits
of the same length then $\alpha (\{ 1, 2, \dots , n-f \})=\{ 1, 2,
\dots , n-f \}$ and $\alpha (\{ n-f+1, \dots , n \})=\{ n-f+1,
\dots , n \}$ and so $\alpha = \gamma \delta$ with $\gamma \in
S\sb{\{ 1, 2, \dots , n-f \}}$ and $\delta \in S\sb{\{ n-f+1,
\dots , n \}}$. In particular, for any $i$ in $\{ 1, 2, \dots ,
n-f \}$, there exists $i'$ in $\{ 1, 2, \dots , n-f \}$ such that
$\alpha (i)=i'$. Then
\begin{align}\notag
(1\; 2\; & \cdots \; n-f)(n-f+1)(n-f+2)\cdots
(n-1)(n)=\mu'\sb{i'}=\mu'\sb{\alpha (i)}=\alpha
\mu\sb{i}\alpha\sp{-1}= \\ \notag &=\gamma\delta (1\; 2\; \cdots
\; n-f)(n-f+1)(n-f+2)\cdots ((n-1)(n)\gamma\sp{-1}\delta\sp{-1}
\end{align}
and so $(1\; 2\; \cdots \; n-f)\gamma =\gamma (1\; 2\; \cdots \;
n-f)$ i.e., $\gamma =(1\; 2\; \cdots \; n-f)\sp{g}$, for some
positive integer $g$. Furthermore, since for any $j$ in $\{ n-f+1,
n-f+2, \dots , n-1, n \}$ there is a $j'$ in $\{ n-f+1, n-f+2,
\dots , n-1, n \}$ such that $\alpha (j)=j'$, then
\begin{align}\notag
(1\; 2\; & \cdots \; n-f)\sp{k'\sb{j'}}(n-f+1)(n-f+2)\cdots
(n-1)(n)=\mu'\sb{j'}=\mu'\sb{\alpha (j)}=\alpha
\mu\sb{j}\alpha\sp{-1}= \\ \notag &=\gamma\delta (1\; 2\; \cdots
\; n-f)\sp{k\sb{j}}(n-f+1)(n-f+2)\cdots
((n-1)(n)\gamma\sp{-1}\delta\sp{-1}
\end{align}
and so
\begin{align}\notag
(1\; 2\; \cdots \; n-f)\sp{k'\sb{\alpha (j)}}&=(1\; 2\; \cdots \;
n-f)\sp{k'\sb{j'}}=(1\; 2\; \cdots \; n-f)\sp{g}(1\; 2\; \cdots \;
n-f)\sp{k\sb{j}}(1\; 2\; \cdots \; n-f)\sp{-g}=\\ \notag &=(1\;
2\; \cdots \; n-f)\sp{k\sb{j}}
\end{align}
i.e., $k'\sb{\alpha (j)}=k\sb{j}$ for any $j$ in $\{ n-f+1, n-f+2,
\dots , n-1, n \}$. Hence, any permutation of the terms of a
sequence of numbers relatively prime to and less than $n-f$, $\{
k\sb{j} \}\sb{j=n-f+1}\sp{n}$, gives rise to isomorphic racks, as
constructed above.

Let us now count how many isomorphism classes there are. Let us
consider them according to how many different numbers there are in
the $\{ k\sb{j} \}\sb{j=n-f+1}\sp{n}$ sequence. If there is just
one number, then there are $\varphi (n-f)$ possibilities, one per
number relatively prime to and less than $n-f$. If there are two
numbers, say $k$ and $k'$, then we may have $f\sb{1}$ repeats of
$k$ and $f\sb{2}$ repeats of $k'$ such that $f\sb{1}+f\sb{2}=f$.
Since we are not allowing permutations of the terms of the
sequence then this is the same as the number of partitions of $f$
using only two summands, $P\sb{2}(f)$. Also the number of
different possibilities of picking two numbers out of a set of
$\varphi (n-f)$ is $C\sb{2}\sp{\varphi (n-f)}$. So the
contribution for the total number of isomorphic classes from using
just two different numbers in the indicated sequence is
$C\sb{2}\sp{\varphi (n-f)}\cdot P\sb{2}(f)$. We believe it is now
clear that when using $i$ distinct numbers the contribution is
$C\sb{i}\sp{\varphi (n-f)}\cdot P\sb{i}(f)$ and that the last $i$
should be the first of $\varphi (n-f)$ or $f$ to be reached. Hence
the total number of isomorphism classes is
$\sum\sb{i=1}\sp{m}C\sb{i}\sp{\varphi (n-f)}P\sb{i}(f)$ where
$m:=\min \{ f, \varphi (n-f) \}$. $\blacksquare$

Proof (Corollary \ref{cor:prackf,n}): Again we assume that the
first $n-f$ permutations are
\[
\mu\sb{1}=\mu\sb{2}=\cdots =\mu\sb{n-f}=(1 \; 2 \; \cdots \;
n-f)(n-f+1)(n-f+2)\dots (n-1)(n)
\]
We now assume one of the remaining permutations fixes the index of
one of the former permutations i.e., there is $j\in \{ n-f+1,
\dots , n \}$ and $i\in \{ 1, \dots , n-f \}$ such that
$\mu\sb{j}(i)=i$. Then
$\mu\sb{i}=\mu\sb{\mu\sb{i}(j)}=\mu\sb{i}\mu\sb{j}\mu\sb{i}\sp{-1}$
i.e., $\mu\sb{i}\mu\sb{j}=\mu\sb{j}\mu\sb{i}$. Hence, mod $n-f$,
we have
\[
\mu\sb{j}(i+1)=\mu\sb{j}\mu\sb{i}(i)=\mu\sb{i}\mu\sb{j}(i)=\mu\sb{i}(i)=i+1
\]
Suppose there is $s$ such that $\mu\sb{j}(i+s)=i+s$, mod $n-f$.
Then, mod $n-f$
\[
\mu\sb{j}(i+s+1)=\mu\sb{j}\mu\sb{i}(i+s)=\mu\sb{i}\mu\sb{j}(i+s)=\mu\sb{i}(i+s)=i+s+1
\]
If the increment is $n-f$ we recover $\mu\sb{j}(i)=i$; then
$\mu\sb{j}$ fixes all elements from $\{ 1, \dots , n-f \}$. If
$f<n-f$ this is impossible. Assume $n-f=f$. Then $\mu\sb{j}$ moves
any element from $\{ n-f+1, \dots , n \}$ so it moves its own
index. Then, without loss of generality, we may assume, invoking
Proposition \ref{prop:perm-deltaorbit}
\[
\mu\sb{n-f+1}=\dots =\mu\sb{n}=(1)\cdots (n-f)(n-f+1 \; n-f+2 \;
\dots \; n-1 \; n)
\]
This prack is decomposable by Proposition \ref{Prop:orbs-decomp}.
It also has two distinct $\Delta$-orbits. Then it is not
isomorphic with any of the pracks in the Theorem since these have
only one $\Delta$-orbit. This finishes the proof.
 $\blacksquare$

\begin{thm}\label{th:prackm1<m2<m3} Any prack of finite order
$m\sb{1}+m\sb{2}+m\sb{3}$ ($1<m\sb{1}<m\sb{2}<m\sb{3}$) with
constant profile
\[
(\{ m\sb{1}, m\sb{2}, m\sb{3}\}, \dots , \{ m\sb{1}, m\sb{2},
m\sb{3}\} )
\]
is decomposable and is isomorphic to a prack given by the
following sequence of permutations,
\begin{align}\notag
\mu\sb{1}&=\dots = \mu\sb{m\sb{1}}=(1\; \cdots \;
m\sb{1})(m\sb{1}+1 \; \cdots \;
m\sb{1}+m\sb{2})\sp{k\sb{2}\sp{1}}(m\sb{1}+m\sb{2}+1 \; \cdots \;
m\sb{1}+m\sb{2}+m\sb{3})\sp{k\sb{3}\sp{1}}\\ \notag
\mu\sb{m\sb{1}+1}&=\cdots = \mu\sb{m\sb{1}+m\sb{2}}=(1\; \cdots \;
m\sb{1})\sp{k\sb{1}\sp{2}}(m\sb{1}+1 \; \cdots \;
m\sb{1}+m\sb{2})(m\sb{1}+m\sb{2}+1 \; \cdots \;
m\sb{1}+m\sb{2}+m\sb{3})\sp{k\sb{3}\sp{2}}\\ \notag
\mu\sb{m\sb{1}+m\sb{2}+1}&=\dots =
\mu\sb{m\sb{1}+m\sb{2}+m\sb{3}}=(1\; \cdots \;
m\sb{1})\sp{k\sb{1}\sp{3}}(m\sb{1}+1 \; \cdots \;
m\sb{1}+m\sb{2})\sp{k\sb{2}\sp{3}}(m\sb{1}+m\sb{2}+1 \; \cdots \;
m\sb{1}+m\sb{2}+m\sb{3})
\end{align}
where $k\sb{1}\sp{2}, k\sb{1}\sp{3}$ are relatively prime to
$m\sb{1}$, $k\sb{2}\sp{1}, k\sb{2}\sp{3}$ are relatively prime to
$m\sb{2}$, and $k\sb{3}\sp{1}, k\sb{3}\sp{2}$ are relatively prime
to $m\sb{3}$.

There are exactly $\bigl( \varphi (m\sb{1})\bigr)\sp{2} \cdot
\bigl( \varphi (m\sb{2})\bigr)\sp{2} \cdot \bigl( \varphi
(m\sb{3})\bigr)\sp{2} $ such pracks up to isomorphism, one
sequence per each $((k\sb{1}\sp{2}, k\sb{1}\sp{3}),
(k\sb{2}\sp{1}, k\sb{2}\sp{3}), (k\sb{3}\sp{1}, k\sb{3}\sp{2}))$,
where $k\sb{1}\sp{2}, k\sb{1}\sp{3}$ are relatively prime to and
less than $m\sb{1}$, $k\sb{2}\sp{1}, k\sb{2}\sp{3}$ are relatively
prime to and less than $m\sb{2}$, and $k\sb{3}\sp{1},
k\sb{3}\sp{2}$ are relatively prime to and less than
$m\sb{3}$.\end{thm}Proof: No element is fixed by any of the
permutations and each of these is the product of three disjoint
cycles, so we can assume, without loss of generality, that
\begin{align}\notag
\mu\sb{1}&=\dots = \mu\sb{m\sb{1}}=(1\; \cdots \; m\sb{1})\pi\sb{2, 3}\\
\notag \mu\sb{m\sb{1}+1}&=\dots =
\mu\sb{m\sb{1}+m\sb{2}}=(m\sb{1}+1 \; \cdots \;
m\sb{1}+m\sb{2})\pi\sb{1, 3}\\ \notag
\mu\sb{m\sb{1}+m\sb{2}+1}&=\dots =
\mu\sb{m\sb{1}+m\sb{2}+m\sb{3}}=(m\sb{1}+m\sb{2}+1 \; \cdots \;
m\sb{1}+m\sb{2}+m\sb{3})\pi\sb{1, 2}
\end{align}
where the $\pi\sb{i, j}$ are from $S\sb{M\sb{i}\cup M\sb{j}}$,
they are products of an $m\sb{i}$-cycle by an $m\sb{j}$-cycle, and
$M\sb{1}=\{ 1, \dots , m\sb{1} \}$, $M\sb{2}=\{ m\sb{1}+1, \dots ,
m\sb{1}+m\sb{2} \}$, and $M\sb{3}=\{ m\sb{1}+m\sb{2}+1, \dots ,
m\sb{1}+m\sb{2}+m\sb{3} \}$.

Let $j$ be in $M\sb{1}$. Since $|M\sb{3}|=m\sb{3}>
m\sb{2}=|M\sb{2}|$, there has to be $i\sb{3}$ in $M\sb {3}$ such
that $\pi\sb{2, 3}(i\sb{3})$ is in $M\sb {3}$. In this way,
\[
(m\sb{1}+m\sb{2}+1 \; \cdots \; m\sb{1}+m\sb{2}+m\sb{3})\pi\sb{1,
2}=\mu\sb{\pi\sb{2,
3}(i\sb{3})}=\mu\sb{\mu\sb{j}(i\sb{3})}=\mu\sb{j}\mu\sb{i\sb{3}}\mu\sb{j}\sp{-1}
\]
i.e.,
\begin{equation}\label{eq:1}
(m\sb{1}+m\sb{2}+1 \; \cdots \; m\sb{1}+m\sb{2}+m\sb{3})\pi\sb{1,
2}=\mu\sb{j}(m\sb{1}+m\sb{2}+1 \; \cdots \;
m\sb{1}+m\sb{2}+m\sb{3})\pi\sb{1, 2}\mu\sb{j}\sp{-1}
\end{equation}

Assume now there is $i\sb{2}$ in $M\sb{2}$ such that $\pi\sb{2,
3}(i\sb{2})$ is in $M\sb{3}$. Then,
\[
(m\sb{1}+m\sb{2}+1 \; \cdots \; m\sb{1}+m\sb{2}+m\sb{3})\pi\sb{1,
2}=\mu\sb{\pi\sb{2,
3}(i\sb{2})}=\mu\sb{\mu\sb{j}(i\sb{2})}=\mu\sb{j}\mu\sb{i\sb{2}}\mu\sb{j}\sp{-1}
\]
i.e.,
\begin{equation}\label{eq:2}
(m\sb{1}+m\sb{2}+1 \; \cdots \; m\sb{1}+m\sb{2}+m\sb{3})\pi\sb{1,
2}=\mu\sb{j}(m\sb{1}+1 \; \cdots \; m\sb{1}+m\sb{2})\pi\sb{1,
3}\mu\sb{j}\sp{-1}
\end{equation}

From (\ref{eq:1}) and (\ref{eq:2}) we have
\[
(m\sb{1}+m\sb{2}+1 \; \cdots \; m\sb{1}+m\sb{2}+m\sb{3})\pi\sb{1,
2}=\mu\sb{i\sb{3}}=\mu\sb{i\sb{2}}=(m\sb{1}+1 \; \cdots \;
m\sb{1}+m\sb{2})\pi\sb{1, 3}
\]
so
\[
\pi\sb{1, 2}=\sigma (m\sb{1}+1 \; \cdots \; m\sb{1}+m\sb{2}),
 \qquad  \qquad \pi\sb{1, 3}=\sigma (m\sb{1}+m\sb{2}+1 \; \cdots \;
m\sb{1}+m\sb{2}+m\sb{3})
\]
with $\sigma$ an $m\sb{1}$-cycle in $S\sb{M\sb{1}}$.

From (\ref{eq:1}) we obtain
\begin{multline*}
(m\sb{1}+m\sb{2}+1 \; \cdots \; m\sb{1}+m\sb{2}+m\sb{3})\cdot
\sigma (m\sb{1}+1 \; \cdots m\sb{1}+m\sb{2})\cdot (1\; \cdots \;
m\sb{1})\pi\sb{2, 3}=\\
=(1\; \cdots \; m\sb{1})\pi\sb{2, 3}\cdot (m\sb{1}+m\sb{2}+1 \;
\cdots \; m\sb{1}+m\sb{2}+m\sb{3})\cdot \sigma (m\sb{1}+1 \;
\cdots \; m\sb{1}+m\sb{2})
\end{multline*}
so $\sigma (1\; \cdots \; m\sb{1})=(1\; \cdots \; m\sb{1})\sigma $
and
\begin{multline}\label{eq:3}
(m\sb{1}+m\sb{2}+1 \; \cdots \; m\sb{1}+m\sb{2}+m\sb{3})(m\sb{1}+1
\; \cdots
\; m\sb{1}+m\sb{2})\pi\sb{2, 3}=\\
=\pi\sb{2, 3} (m\sb{1}+m\sb{2}+1 \; \cdots \;
m\sb{1}+m\sb{2}+m\sb{3})(m\sb{1}+1 \; \cdots \; m\sb{1}+m\sb{2})
\end{multline}

We recall we are assuming there exists $i\sb{2}$ in $M\sb{2}$ such
that $\pi\sb{2, 3}(i\sb{2})$ is in $M\sb{3}$. Then, using
(\ref{eq:3}), we have mod $m\sb{3}$,
\begin{multline*}
\pi\sb{2, 3}(i\sb{2}+1)=\pi\sb{2, 3} (m\sb{1}+m\sb{2}+1 \; \cdots
\; m\sb{1}+m\sb{2}+m\sb{3})(m\sb{1}+1 \; \cdots
\; m\sb{1}+m\sb{2})(i\sb{2})=\\
=(m\sb{1}+m\sb{2}+1 \; \cdots \;
m\sb{1}+m\sb{2}+m\sb{3})(m\sb{1}+1 \; \cdots \;
m\sb{1}+m\sb{2})\pi\sb{2, 3}(i\sb{2})=\pi\sb{2, 3}(i\sb{2})+1
\end{multline*}
Now suppose there is $s$ such that $\pi\sb{2,
3}(i\sb{2}+s)=\pi\sb{2, 3}(i\sb{2})+s$ in $M\sb{3}$. Then,
\begin{multline*}
\pi\sb{2, 3}(i\sb{2}+s+1)=\pi\sb{2, 3} (m\sb{1}+m\sb{2}+1 \;
\cdots \; m\sb{1}+m\sb{2}+m\sb{3})(m\sb{1}+1 \; \cdots
\; m\sb{1}+m\sb{2})(i\sb{2}+s)=\\
=(m\sb{1}+m\sb{2}+1 \; \cdots \;
m\sb{1}+m\sb{2}+m\sb{3})(m\sb{1}+1 \; \cdots \;
m\sb{1}+m\sb{2})\pi\sb{2, 3}(i\sb{2}+s)=\pi\sb{2,
3}(i\sb{2}+s)+1=\pi\sb{2, 3}(i\sb{2})+s+1
\end{multline*}
In particular,
\[
\pi\sb{2, 3}(i\sb{2})=\pi\sb{2, 3}(i\sb{2}+m\sb{2})=\pi\sb{2,
3}(i\sb{2})+m\sb{2}
\]
so $m\sb{2}=m\sb{3}$ mod $m\sb{3}$ but by hypothesis,
$1<m\sb{2}<m\sb{3}$. Hence, there cannot be $i\sb{2}$ in $M\sb{2}$
such that $\pi\sb{2, 3}(i\sb{2})$ is in $M\sb{3}$. So
\[
\pi\sb{2, 3}=\gamma\sb{2}\sp{1} \delta\sb{3}\sp{1}
\]
where $\gamma\sb{2}\sp{1}$ is an $m\sb{2}$-cycle in
$S\sb{M\sb{2}}$, and $\delta\sb{3}\sp{1} $ is an $m\sb{3}$-cycle
in $S\sb{M\sb{3}}$.

With similar arguments we obtain also that
\[
\pi\sb{1, 3}=\gamma\sb{1}\sp{2} \delta\sb{3}\sp{2}
\]
where $\gamma\sb{1}\sp{2} $ is an $m\sb{1}$-cycle in
$S\sb{M\sb{1}}$, and $\delta\sb{3}\sp{2} $ is an $m\sb{3}$-cycle
in $S\sb{M\sb{3}}$, and
\[
\pi\sb{1, 2}=\gamma\sb{1}\sp{3} \delta\sb{2}\sp{3}
\]
where $\gamma\sb{1}\sp{3} $ is an $m\sb{1}$-cycle in
$S\sb{M\sb{1}}$, and $\delta\sb{2}\sp{3}$ is an $m\sb{2}$-cycle in
$S\sb{M\sb{2}}$.

So, for any $j\sb{1}$ in $M\sb{1}$, $j\sb{2}$ in $M\sb{2}$ and
$j\sb{3}$ in $M\sb{3}$, $\mu\sb{j\sb{1}}(j\sb{2})$ is in $M\sb{2}$
and $\mu\sb{j\sb{1}}(j\sb{3})$ is in $M\sb{3}$, so
\begin{align}\notag
(m\sb{1}+1\; \cdots \;
m\sb{1}+m\sb{2})\gamma\sb{1}\sp{2} \delta\sb{3}\sp{2}&=\mu\sb{j\sb{1}}(j\sb{2})=\mu\sb{j\sb{1}}\mu\sb{j\sb{2}}\mu\sb{j\sb{1}}\sp{-1}\\
\notag &=(1\; \cdots \;
m\sb{1})\gamma\sb{2}\sp{1}\delta\sb{3}\sp{1}\cdot (m\sb{1}+1\;
\cdots \; m\sb{1}+m\sb{2})\gamma\sb{1}\sp{2}
\delta\sb{3}\sp{2}\cdot (1\; \cdots \;
m\sb{1})\sp{-1}{\gamma\sb{2}\sp{1}}sp{-1}{\delta\sb{3}\sp{1}}\sp{-1}
\end{align}
i.e.,
\[
(m\sb{1}+1\; \cdots \; m\sb{1}+m\sb{2})\gamma\sb{1}\sp{2}
\delta\sb{3}\sp{2}(1\; \cdots \;
m\sb{1})\gamma\sb{2}\sp{1}\delta\sb{3}\sp{1} =(1\; \cdots \;
m\sb{1})\gamma\sb{2}\sp{1}\delta\sb{3}\sp{1}(m\sb{1}+1\; \cdots \;
m\sb{1}+m\sb{2})\gamma\sb{1}\sp{2}  \delta\sb{3}\sp{2}
\]
so
\[
\gamma\sb{1}\sp{2}  (1\; \cdots \; m\sb{1})=(1\; \cdots \;
m\sb{1})\gamma\sb{1}\sp{2}   , \qquad
\gamma\sb{2}\sp{1}(m\sb{1}+1\; \cdots \;
m\sb{1}+m\sb{2})=(m\sb{1}+1\; \cdots \;
m\sb{1}+m\sb{2})\gamma\sb{2}\sp{1}, \qquad
\delta\sb{3}\sp{1}\delta\sb{3}\sp{2}=\delta\sb{3}\sp{2}\delta\sb{3}\sp{1}
\]
So, in particular, $\gamma\sb{1}\sp{2}  $ is an $m\sb{1}$-cycle in
the centralizer of $(1\; \cdots \; m\sb{1})$ in $S\sb{M\sb{1}}$
i.e., $\gamma\sb{1}\sp{2}  =(1\; \cdots \;
m\sb{1})\sp{k\sb{1}\sp{2}}$, where $k\sb{1}\sp{2}$ is relatively
prime to $m\sb{1}$. Also, $\gamma\sb{2}\sp{1}$ is an
$m\sb{2}$-cycle in the centralizer of $(m\sb{1}+1\; \cdots \;
m\sb{1}+m\sb{2})$ in $S\sb{M\sb{2}}$, thus
$\gamma\sb{2}\sp{1}=(m\sb{1}+1\; \cdots \;
m\sb{1}+m\sb{2})\sp{k\sb{2}\sp{1}}$, with $k\sb{2}\sp{1}$
relatively prime to $m\sb{2}$. Repeating this argument with other
pairs of indices we eventually obtain $\delta\sb{3}\sp{1} =
(m\sb{1}+m\sb{2}+1\; \cdots \;
m\sb{1}+m\sb{2}+m\sb{3})\sp{k\sb{3}\sp{1}}$, with $k\sb{3}\sp{1}$
relatively prime to $m\sb{3}$, $\delta\sb{3}\sp{2} =
(m\sb{1}+m\sb{2}+1\; \cdots \;
m\sb{1}+m\sb{2}+m\sb{3})\sp{k\sb{3}\sp{2}}$, with $k\sb{3}\sp{2}$
relatively prime to $m\sb{3}$, $\gamma\sb{1}\sp{3}  = (1\; \cdots
\; m\sb{1})\sp{k\sb{1}\sp{3}}$, with $k\sb{1}\sp{3}$ relatively
prime to $m\sb{1}$, and $\delta\sb{2}\sp{3} = (m\sb{1}+1\; \cdots
\; m\sb{1}+m\sb{2})\sp{k\sb{2}\sp{3}}$, with $k\sb{2}\sp{3}$
relatively prime to $m\sb{2}$ which gives the permutations in the
statement. In particular, such a rack is decomposable since the
elements which form the cycles the permutations are made of split
into two disjoint subsets, see Proposition \ref{Prop:orbs-decomp}.

We remark that no two of these pracks are isomorphic, when the
$k\sb{i}\sp{j}$'s are picked less than the corresponding
$m\sb{i}$'s. The proof goes along the same lines as in the proofs
to previous theorems. An isomorphism between two such pracks maps
$\Delta$-orbits into $\Delta$-orbits, so it will be the product of
three permutations each of which will map $M\sb{i}$ into itself,
for each $i$. Next each of these permutation has to be in a
certain centralizer. Finally this will imply that the
corresponding $k\sb{i}\sp{j}$ exponents have to be equal which
means the pracks are equal. It then follows that the number of
distinct isomorphic classes is $\bigl( \varphi
(m\sb{1})\bigr)\sp{2} \cdot \bigl( \varphi (m\sb{2})\bigr)\sp{2}
\cdot \bigl( \varphi (m\sb{3})\bigr)\sp{2} $, which concludes the
proof.
 $\blacksquare$

We turn now to quandles of specified profile. Clearly, there
cannot be quandles of constant profile $(\{ n \},\dots ,\{ n \})$.
What about quandles of constant profile $(\{1, n-1\},\dots , \{1,
n-1\})$?

\begin{thm}\label{th:1,n-1} Consider a finite quandle of order $n>2$ with constant profile
$(\{1, n-1\},\dots , \{1, n-1\})$. Then, modulo isomorphism, its
sequence of permutations is:
\begin{enumerate}
\item  $\mu\sb{n}=( 1\; 2 \; \cdots \; n-1 )(n)$,

\item $\mu\sb{n-1}$ is a solution to the system of equations:
\begin{itemize}
\item
$\mu\sb{n-1}\mu\sb{n}\mu\sb{n-1}\sp{-1}=\mu\sb{n}\sp{\mu\sb{n-1}(n)}\mu\sb{n-1}\mu\sb{n}\sp{-\mu\sb{n-1}(n)}$
\item
$\mu\sb{n-1}\sp{-1}\mu\sb{n}\mu\sb{n-1}=\mu\sb{n}\sp{\mu\sp{-1}\sb{n-1}(n)}\mu\sb{n-1}\mu\sb{n}\sp{-\mu\sp{-1}\sb{n-1}(n)}$
\item
$\mu\sb{n}\sp{-\mu\sb{n-1}(l)}\mu\sb{n-1}\mu\sb{n}\sp{l}=\mu\sb{n-1}\sp{k\sb{l}}$,
for all $l \in \{ 1, \dots , n-2 \}\setminus \{
\mu\sb{n-1}\sp{-1}(n) \}$, some $k\sb{l} \in \{ 1,\dots , n-2 \}$.
\end{itemize}
\item
$\mu\sb{k}=\mu\sb{n}\sp{k}\mu\sb{n-1}\mu\sb{n}\sp{-k}$, for all
$1\leq k \leq n-2$
\end{enumerate}
\end{thm}Proof: Without loss of generality, assume $\mu\sb{n}=(1 \; 2 \; \cdots \;
n-1)(n)$, whence $1.$ follows. Now, $\mu\sb{n-1}=(i\sb{1}\;
i\sb{2} \; \cdots \; i\sb{n-1})(n-1)$, where $(i\sb{1}\; i\sb{2}
\; \cdots \; i\sb{n-1})$ is a cycle of length $n-1$, permuting the
elements of $\{ 1, 2,\dots , n-2, n-1, n \}\setminus \{ n-1 \}$.
Using (part of) the expressions corresponding to the
self-distributivity axiom, we can define the remaining
permutations. We have
\[
\mu\sb{1}=\mu\sb{\mu\sb{n}(n-1)}=\mu\sb{n}\mu\sb{n-1}\mu\sb{n}\sp{-1}
\]
If $\mu\sb{k}=\mu\sb{n}\sp{k}\mu\sb{n-1}\mu\sb{n}\sp{-k}$, then
\[
\mu\sb{k+1}=\mu\sb{\mu\sb{n}(k)}=\mu\sb{n}\mu\sb{k}\mu\sb{n}\sp{-1}=\mu\sb{n}\mu\sb{n}\sp{k}\mu\sb{n-1}\mu\sb{n}\sp{-k}\mu\sb{n}\sp{-1}=\mu\sb{n}\sp{k+1}\mu\sb{n-1}\mu\sb{n}\sp{-(k+1)}
\]
(where the free indices are to be read modulo $n-1$) whence $3.$
follows. Note that for $1\leq k \leq n-2$ we have
\begin{align}\notag
\mu\sb{k}(k)&=\mu\sb{n}\sp{k}\mu\sb{n-1}\mu\sb{n}\sp{-k}(k)=\mu\sb{n}\sp{k}\mu\sb{n-1}\mu\sb{n}\sp{-k+1}(k-1)=\dots
=\mu\sb{n}\sp{k}\mu\sb{n-1}\mu\sb{n}\sp{-1}(1)=\mu\sb{n}\sp{k}\mu\sb{n-1}(n-1)=\\
\notag &=\mu\sb{n}\sp{k}(n-1)=\mu\sb{n}\sp{k-1}(1)=\dots
=\mu\sb{n}(k-1)=k
\end{align}
so these permutations comply with the quandle axiom. Now for the
self-distributivity,
$\mu\sb{\mu\sb{i}(j)}=\mu\sb{i}\mu\sb{j}\mu\sb{i}\sp{-1}$. Before,
we used these expressions for $i=n$, and recursively $j=n-1, 1,
2,\dots n-3$, for defining $\mu\sb{k}$, $1\leq k \leq n-2$. For
$k$ and $k'$ in $\{ 1, 2,\dots , n-1 \}$ ($k\neq k'$) we must have
\[
\mu\sb{\mu\sb{k}(k')}=\mu\sb{k}\mu\sb{k'}\mu\sb{k}\sp{-1}
\]
Now,
\[
\mu\sb{\mu\sb{k}(k')}=\mu\sb{\mu\sb{n}\sp{k}\mu\sb{n-1}\mu\sb{n}\sp{-k}(k')}=\mu\sb{\mu\sb{n}\sp{k}\mu\sb{n-1}(k'-k)}=\mu\sb{\mu\sb{n}\sp{k}(\mu\sb{n-1}(l))}=
\]
\begin{equation}\notag
=
\begin{cases}
\mu\sb{n}, &\text{ if $\mu\sb{n-1}(l)=n$}\\ \notag
\mu\sb{n}\sp{k}\mu\sb{\mu\sb{n-1}(l)}\mu\sb{n}\sp{-k}, &\text{ if
$\mu\sb{n-1}(l)\in \{ 1,\dots , n-2 \}$}
\end{cases}
\end{equation}
with $l=k'-k$. On the other hand, again with $l=k'-k$,
\[
\mu\sb{k}\mu\sb{k'}\mu\sb{k}\sp{-1}=\bigl(
\mu\sb{n}\sp{k}\mu\sb{n-1}\mu\sb{n}\sp{-k}\bigr)\cdot \bigl(
\mu\sb{n}\sp{k'}\mu\sb{n-1}\mu\sb{n}\sp{-k'}\bigr)\cdot \bigl(
\mu\sb{n}\sp{k}\mu\sb{n-1}\mu\sb{n}\sp{-k}\bigr)\sp{-1}=\mu\sb{n}\sp{k}\mu\sb{n-1}\mu\sb{n}\sp{l}\mu\sb{n-1}\mu\sb{n}\sp{-l}
\mu\sb{n-1}\sp{-1}\mu\sb{n}\sp{-k}
\]

Suppose $\mu\sb{n-1}(l)=n$. Then,
$\mu\sb{n}=\mu\sb{n-1}\mu\sb{n}\sp{l}\mu\sb{n-1}\mu\sb{n}\sp{-l}\mu\sb{n-1}\sp{-1}$,
i.e.,
\begin{equation}\label{eq:2.2}
\mu\sb{n-1}\sp{-1}\mu\sb{n}\mu\sb{n-1}=\mu\sb{n}\sp{l}\mu\sb{n-1}\mu\sb{n}\sp{-l},\qquad
 \text{ if $\mu\sb{n-1}(l)=n$ }
\end{equation}

Suppose $\mu\sb{n-1}(l) \in \{ 1, 2, \dots , n-2 \}$. Then,
$\mu\sb{\mu\sb{n-1}(l)}=\mu\sb{n-1}\mu\sb{n}\sp{l}\mu\sb{n-1}\mu\sb{n}\sp{-l}\mu\sb{n-1}\sp{-1}$.
On the other hand, since $\mu\sb{n-1}(l) \in \{ 1, 2, \dots , n-2
\}$,
$\mu\sb{\mu\sb{n-1}(l)}=\mu\sb{n}\sp{\mu\sb{n-1}(l)}\mu\sb{n-1}\mu\sb{n}\sp{-\mu\sb{n-1}(l)}
$ and so
\[
\mu\sb{n}\sp{\mu\sb{n-1}(l)}\mu\sb{n-1}\mu\sb{n}\sp{-\mu\sb{n-1}(l)}=\mu\sb{n-1}\mu\sb{n}\sp{l}\mu\sb{n-1}\mu\sb{n}\sp{-l}\mu\sb{n-1}\sp{-1}
\]
i.e.,
\[
\mu\sb{n-1}\cdot
\mu\sb{n}\sp{-\mu\sb{n-1}(l)}\mu\sb{n-1}\mu\sb{n}\sp{l}=\mu\sb{n}\sp{-\mu\sb{n-1}(l)}\mu\sb{n-1}\mu\sb{n}\sp{l}\cdot
\mu\sb{n-1}
\]
This means
$\mu\sb{n}\sp{-\mu\sb{n-1}\sp{l}}\mu\sb{n-1}\mu\sb{n}\sp{l}$ is in
the centralizer of $\mu\sb{n-1}$ in the symmetry group of $\{ 1,
2, \dots , n \}$ which is the cyclic subgroup generated by
$\mu\sb{n-1}$. Thus,
\begin{equation}\label{eq:2.3}
\mu\sb{n}\sp{-\mu\sb{n-1}(l)}\mu\sb{n-1}\mu\sb{n}\sp{l}=\mu\sb{n-1}\sp{k\sb{l}},\quad
 \text{ for all $l \in  \{1, \dots , n-2\}\setminus \{ \mu\sb{n-1}\sp{-1}(n) \}$, and some $k\sb{l} \in  \{1, \dots , n-2\}$  }
\end{equation}
We remark that $k\sb{l}=n-1$ leads to a contradiction.

There remains to ascertain the implications of
$\mu\sb{\mu\sb{k}(n)}=\mu\sb{k}\mu\sb{n}\mu\sb{k}\sp{-1}$, for $k$
in $ \{1, \dots , n-1\}$. With similar manipulations as above,
noting that
\[
\mu\sb{k}(n)=\mu\sb{n}\sp{k}\mu\sb{n-1}\mu\sb{n}\sp{-k}(n)=\mu\sb{n}\sp{k}\mu\sb{n-1}(n)
=\mu\sb{n}\sp{k}\mu\sb{n}\sp{\mu\sb{n-1}(n)}(n-1)=\mu\sb{n}\sp{k+\mu\sb{n-1}(n)}(n-1)
\]
we obtain,
\begin{equation}\label{eq:2.1}
\mu\sb{n-1}\mu\sb{n}\mu\sb{n-1}\sp{-1}=\mu\sb{n}\sp{\mu\sb{n-1}(n)}\mu\sb{n-1}\mu\sb{n}\sp{-\mu\sb{n-1}(n)}
\end{equation}
(\ref{eq:2.1}), (\ref{eq:2.2}), and (\ref{eq:2.3}) now stand for
statement $2.$. This concludes the proof.
 $\blacksquare$

The following Corollaries depict restrictions the $\mu\sb{n-1}$ is
subject to. In these Corollaries the context and notation are
those of Theorem \ref{th:1,n-1}.

\begin{cor} $\mu\sb{n-1}(i)\neq i+1$, for $i$ in $\{ 1, 2,\dots , n-2 \}$.
\end{cor}Proof: Assume it is. Then
$\mu\sb{i+1}=\mu\sb{\mu\sb{n-1}(i)}=\mu\sb{n-1}\mu\sb{i}\mu\sb{n-1}\sp{-1}$.
Since $\mu\sb{i}=\mu\sb{n}\sp{i}\mu\sb{n-1}\mu\sb{n}\sp{-i}$ and
analogously for $\mu\sb{i+1}$, then
$\mu\sb{n}\sp{i+1}\mu\sb{n-1}\mu\sb{n}\sp{-(i+1)}=\mu\sb{n-1}\mu\sb{n}\sp{i}\mu\sb{n-1}\mu\sb{n}\sp{-i}\mu\sb{n-1}\sp{-1}$
which can be rewritten as $\mu\sb{n-1}\mu\sb{n}\sp{-1}\cdot
\mu\sb{n}\sp{-i}\mu\sb{n-1}\mu\sb{n}\sp{i}=\mu\sb{n}\sp{-1}\cdot
\mu\sb{n}\sp{-i}\mu\sb{n-1}\mu\sb{n}\sp{i}\cdot \mu\sb{n-1}$,
which means
$\mu\sb{n-1}\mu\sb{n}\sp{-1}\mu\sb{-i}=\mu\sb{n}\sp{-1}\mu\sb{-i}\mu\sb{n-1}$,
($-i$ mod $n-1$) i.e., $\mu\sb{n}\sp{-1}\mu\sb{-i}$ is in the
centralizer of $\mu\sb{n-1}$ in the symmetry group of $\{ 1, 2,
\dots , n \}$. Hence,
$\mu\sb{n}\sp{-1}\mu\sb{-i}=\mu\sb{n-1}\sp{k}$, for some $k$ in
$\{ 1, 2, \dots , n-1 \}$. If $k=n-1$, then $\mu\sb
{-i}=\mu\sb{n}$. Thus, modulo $n-1$, $-i=n$ which means that
$i=n-2$. So $\mu\sb{n-1}(n-2)=n-1$ but by hypothesis,
$\mu\sb{n-1}$ fixes $n-1$. This is impossible.

Assume now $k\neq n-1$. Then, since $\mu\sb{n-1}\sp{k}$ fixes
$n-1$ and $\mu\sb{n}\sp{-1}\mu\sb{-i}$ fixes
$\mu\sb{-i}\sp{-1}(n)$, hence $\mu\sb{-i}\sp{-1}(n)=n-1$, i.e.,
$n=\mu\sb{-i}(n-1)=\mu\sb{n}\sp{-i}\mu\sb{n-1}\mu\sb{n}\sp{i}(n-1)=\mu\sb{n}\sp{-i}\mu\sb{n-1}(i)$
from which follows that $n=\mu\sb{n-1}(i)=i+1$, by assumption, so
$n-1 = i \in \{ 1, 2,\dots , n-2 \}$. Again this is impossible.

 $\blacksquare$

We remark that $\mu\sb{n-1}(n-2)=n$ is possible: consider
$\mathbb{Z}\sb{5}[T, T\sp{-1}]/(T-2)$; and that $\mu\sb{n-1}(n)=1$
is also possible: $\mathbb{Z}\sb{5}[T, T\sp{-1}]/(T-3)$.

\begin{cor} If, for $n>3$,
\[
\mu\sb{n-1}=( n-2 \; n \; k \; \cdots )( n-1 )
\]
then $n$ is odd and $k=\frac{n-1}{2}$.
\end{cor}Proof: Suppose $\mu\sb{n-1}$ is as stated. Then,
\[
\mu\sb{n-1}\mu\sb{n}\mu\sb{n-1}\sp{-1}=\bigl( \mu\sb{n-1}(1) \;
 \mu\sb{n-1}(2) \;  \cdots \;  \mu\sb{n-1}(n-2) \;
 \mu\sb{n-1}(n-1) \;
 \;
  \bigr)  \bigl( \mu\sb{n-1}(n)  \bigr) = \bigl( \; \cdots
\;  n \; n-1 \;  \bigr)  \bigl( k  \bigr)
\]
On the other hand,
\[
\mu\sb{n}\sp{\mu\sb{n-1}(n)}\mu\sb{n-1}\mu\sb{n}\sp{-\mu\sb{n-1}(n)}=\mu\sb{n}\sp{k}\mu\sb{n-1}\mu\sb{n}\sp{-k}=\bigl(
\mu\sb{n}\sp{k}(n-2) \;
 \mu\sb{n}\sp{k}(n) \; \mu\sb{n}\sp{k}(k) \; \cdots  \bigr)  \bigl( \mu\sb{n}\sp{k}(n-1)  \bigr) = \bigl( k-1 \;
 n \; 2k \;  \cdots
  \bigr)  \bigl( k  \bigr)
\]
By (\ref{eq:2.1}) in the proof to the previous Theorem, these two
expressions have to be equal. Then, we have $2k=n-1$ i.e., $n$ is
odd and $k=(n-1)/2$

$\blacksquare$

\begin{cor}\label{cor:s4} If, for $n>3$,
\[
\mu\sb{n-1}=( k\;  n \; n-2 \; \cdots )( n-1 )
\]
 then $k=1$.
\end{cor}Proof: Assume $\mu\sb{n-1}$ is as stated. Then
\[
\mu\sb{n-1}\sp{-1}\mu\sb{n}\mu\sb{n-1}=\bigl(
\mu\sb{n-1}\sp{-1}(1) \;
 \cdots  \; \mu\sb{n-1}\sp{-1}(n-2) \; \cdots \;
 \mu\sb{n-1}\sp{-1}(n-1) \;
 \;
  \bigr)  \bigl( \mu\sb{n-1}\sp{-1}(n)  \bigr) = \bigl(  \;  \cdots
  \;
 n \; n-1 \; \bigr)  \bigl( k  \bigr)
\]
On the other hand,
\begin{align}\notag
\mu\sb{n}\sp{\mu\sb{n-1}\sp{-1}(n)}\mu\sb{n-1}\mu\sb{n}\sp{-\mu\sb{n-1}\sp{-1}(n)}&=\mu\sb{n}\sp{k}\mu\sb{n-1}\mu\sb{n}\sp{-k}=\bigl(
\mu\sb{n}\sp{k}(k) \;
 \mu\sb{n}\sp{k}(n) \; \mu\sb{n}\sp{k}(n-2) \; \cdots  \bigr)  \bigl( \mu\sb{n}\sp{k}(n-1)  \bigr) = \\ \notag
 &=\bigl( 2k \;
 n \; k-1 \;  \cdots
  \bigr)  \bigl( k  \bigr)
\end{align}
and by comparison, $k=1$.
 $\blacksquare$

The proofs of the following two Corollaries are similar to the
ones in the preceding corollaries and so will be omitted.

\begin{cor} If, for $n>3$,
\[
\mu\sb{n-1}=( k\;  n \; 1 \; \cdots )( n-1 )
\]
 then $n$ is odd and $k=\frac{n-1}{2}$.
\end{cor}
 $\blacksquare$

\begin{cor} If, for $n>3$,
\[
\mu\sb{n-1}=( 1\;  n \; k \; \cdots )( n-1 )
\]
 then $k=n-2$.
\end{cor}
 $\blacksquare$

Using these results we tried to find $\mu\sb{n-1}$'s which are
solutions to the system of equations in $2.$ in the statement of
Theorem \ref{th:1,n-1}. For each $n>2$, each of these solutions
together with $\mu\sb{n}=( 1\; 2\dots n-1 )(n)$ and the
$\mu\sb{k}$'s given by
$\mu\sb{k}=\mu\sb{n}\sp{k}\mu\sb{n-1}\mu\sb{n}\sp{-k}$, for all
$1\leq k \leq n-2$, constitute the sequence of permutations of a
finite quandle of order $n>2$ with constant profile $(\{1,
n-1\},\dots , \{1, n-1\})$. The general strategy was, for each
order $n$, to set up a tentative $\mu\sb{n-1}$ by assigning the
image and pre-image of $n$, taking into consideration the
Corollaries above. Then, using the first two equations of the
system of equations referred to above we tried either to show that
such a $\mu\sb{n-1}$ would not be a solution to our problem or to
complete $\mu\sb{n-1}$ and to verify that it satisfied the
remaining equations.

For $n=3$, the unique such $\mu\sb{n-1}$ is $\mu\sb{2}=(1\; 3)( 2
)$ which gives rise to the dihedral quandle of order $3$,
$R\sb{3}$. Applying Corollary \ref{cor:s4}, with $n=4$, the unique
such $\mu\sb{n-1}$ is $\mu\sb{3}=( 1\; 4 \; 2\; )(3)$ which gives
rise to a quandle isomorphic to $S\sb{4} \cong \mathbb{Z}\sb{2}[T,
T\sp{-1}]/(T\sp{2}+T+1)$. And so on and so forth. We display our
results in Table \ref{Ta:1,n-1}. The order is displayed in the
first column. In the second column the $\mu\sb{n-1}$ is displayed
since it identifies the quandle under study (see discussion
above). In the third column, we give a more familiar quandle which
is isomorphic with the quandle under study.

\begin{table}
\begin{center}
\renewcommand{\arraystretch}{1.25}
    \begin{tabular}{| c | c | c |}\hline
Order          &        $\mu\sb{n-1}$          & ...isomorphic
with...             \\ \hline
  $3$     &       $(2)(1 \; 3)$             & $R\sb{3}$            \\ \hline
  $4$     &      $( 3 )( 1 \; 4 \; 2 )$     & $S\sb{4}\cong \mathbb{Z}\sb{2}[T, T\sp{-1}]/(T\sp{2}+T+1)$      \\ \hline
  $5$     &      $( 4 )( 2 \; 5 \; 1 \; 3 )$&
$\mathbb{Z}\sb{5}[T, T\sp{-1}]/(T-3)$   \\  \cline{2-3}
  &  $( 4 )( 3 \; 5 \; 2 \; 1 )$      &     $\mathbb{Z}\sb{5}[T, T\sp{-1}]/(T-2)$            \\\hline
 $6$    & no solutions    &                     \\ \hline
 $7$     & $( 6 )( 1 \; 7 \; 5 \; 2 \; 4 \; 3 )$ &
$\mathbb{Z}\sb{7}[T, T\sp{-1}]/(T-5)$\\ \cline{2-3}
 &     $( 6 )( 1 \; 7 \; 5 \; 3  \; 2 \; 4 )$     &                      $\mathbb{Z}\sb{7}[T, T\sp{-1}]/(T-3)$      \\\hline
 $8$     & $( 7 )( 2 \; 8 \; 3 \; 6 \; 1 \; 5 \; 4 )$ &
$\mathbb{Z}\sb{2}[T, T\sp{-1}]/(T\sp{3}+T+1)$\\ \cline{2-3}
 &     $( 7 )( 4 \; 8 \; 5 \; 3  \; 2 \; 6 \; 1 )$     &                      $\mathbb{Z}\sb{2}[T, T\sp{-1}]/(T\sp{3}+T\sp{2}+1)$      \\\hline
    \end{tabular}
    \caption{Quandles of profile $(\{ 1, n-1 \}, \dots , \{ 1, n-1 \})$ for each order $3\leq n\leq 8$}\label{Ta:1,n-1}
\end{center}
\end{table}

\section{Final Remarks} \label{S:final}

In this work we developed a different approach to racks and
quandles by regarding them as sequences of permutations. In the
course of it, some issues concerning the structure of racks and
quandles were found. In Proposition \ref{prop: aqsd} it is proved
that the opposite to an Alexander quandle also satisfies
self-distributivity. We feel this is an indication to look for
non-Alexander quandles. We have a first example in the octahedron
quandle. This is an interesting quandle in that it is good at
telling knots apart via counting colorings (\cite{D}). This is
also an indecomposable quandle. In this way, we would like to know
if there is an infinite class of indecomposable non-Alexander
quandles which are good at telling knots apart via counting
colorings. Finally, we were able to identify all racks that have
given profiles. In this connection we ask the following question.
Are there indecomposable racks not isomorphic with cyclic racks?

\subsection{Acknowledgements} \label{Ss:ack}
The first author acknowledges support by {\em Programa Operacional
``Ci\^{e}ncia, Tecnologia, Inova\c{c}\~{a}o''} (POCTI) of the {\em
Funda\c{c}\~{a}o para a Ci\^{e}ncia e a Tecnologia} (FCT)
cofinanced by the European Community fund FEDER.

\end{document}